\documentclass[]{siamart}

\usepackage[capitalize]{cleveref}

\crefformat{section}{#2\S#1#3}
\crefrangeformat{section}{#3\S\S#1#4--#5#2#6}
\crefmultiformat{section}{#2\S\S#1#3}{ and #2#1#3}{, #2#1#3}{, and #2#1#3}
\crefrangemultiformat{section}{#3\S\S#1#4--#5#2#6}%
{ and #3#1#4--#5#2#6}{, #3#1#4--#5#2#6}{, and #3#1#4--#5#2#6}

\Crefformat{section}{#2Section~#1#3}
\Crefrangeformat{section}{#3Sections~#1#4--#5#2#6}
\Crefmultiformat{section}{#2Sections~#1#3}{ and #2#1#3}{, #2#1#3}{, and #2#1#3}
\Crefrangemultiformat{section}{#3Sections~#1#4--#5#2#6}%
{ and #3#1#4--#5#2#6}{, #3#1#4--#5#2#6}{, and #3#1#4--#5#2#6}

\crefformat{equation}{#2(#1)#3}
\crefrangeformat{equation}{(#3#1#4--#5#2#6)}
\crefmultiformat{equation}{#2(#1)#3}{ and #2(#1)#3}
{, #2(#1)#3}{, and #2(#1)#3}
\crefrangemultiformat{equation}{(#3#1#4--#5#2#6)}%
{ and (#3#1#4--#5#2#6)}{, (#3#1#4--#5#2#6)}{, and (#3#1#4--#5#2#6)}

\Crefformat{equation}{#2Equation~(#1)#3}
\Crefrangeformat{equation}{Equations~(#3#1#4--#5#2#6)}
\Crefmultiformat{equation}{Equations~(#2#1#3)}{ and (#2#1#3)}
{, (#2#1#3)}{, and (#2#1#3)}
\Crefrangemultiformat{equation}{Equations~(#3#1#4--#5#2#6)}%
{ and (#3#1#4--#5#2#6)}{, (#3#1#4--#5#2#6)}{, and (#3#1#4--#5#2#6)}

\usepackage{color}
\definecolor{red}{rgb}{1.0,0.,0.}

\usepackage{amsmath}
\usepackage{amssymb}
\usepackage{multicol}

\usepackage{verbatim}
\usepackage{accents}
\usepackage{graphicx}
\usepackage{subfigure}
\usepackage{multirow}

\newcommand\reals{{{\rm l} \kern -.15em {\rm R} }}
\newcommand{\ignore}[1]{}

\graphicspath{.}

\begin{document}

\title{Adjoint Methods for Guiding Adaptive Mesh Refinement in Wave
Propagation Problems%
  \thanks{Supported in part by an NSF Graduate Research Fellowship 
  DGE-1256082 and NSF grants DMS-1216732 and EAR-133141.}}

\subtitle{\today} 

\author{Brisa N. Davis
  \thanks{Department of Applied Mathematics, University of Washington,
Seattle, WA \email{bndavis@uw.edu}.}%
  \and
  Randall J. LeVeque%
  \thanks{Department of Applied Mathematics, University of Washington,
Seattle, WA \email{rjl@uw.edu}.}%
}

\maketitle


\pagestyle{myheadings}
\thispagestyle{plain}
\markboth{B. N. Davis and R. J. LeVeque}%
{Adjoint Guided Adaptive Mesh Refinement}

\begin{abstract}
One difficulty in developing numerical methods for
hyperbolic systems of conservation laws is the fact 
that solutions often contain regions where much higher resolution is required
than elsewhere in the domain, particularly since the solution may contain
discontinuities or other localized features.
The Clawpack software deals with 
this issue by using block-structured adaptive mesh 
refinement to selectively refine around propagating waves. For problems 
where only a target area of the total solution is of interest, a method that 
allows identifying and refining the grid only in regions that influence this 
target area would significantly reduce the computational cost of finding a 
solution. 

In this work, we show that solving the time-dependent adjoint equation and using 
a suitable inner product with the forward solution allows more precise 
refinement of the relevant waves. We present acoustics examples in one and two
dimensions and a tsunami propagation example. To perform these simulations,
the use of the adjoint 
method has been integrated into the adaptive mesh refinement strategy of
the open source Clawpack and GeoClaw software.
We also present results that show that the accuracy of the solution is maintained 
and the computational time required is significantly reduced through the integration 
of the adjoint method into AMR.
\end{abstract}



\section{Introduction}
\label{sec:intro}

\indent Hyperbolic systems of conservation laws 
\begin{equation}\label{eq:hypsys}
q_t + f(q)_x + g(q)_y = 0
\end{equation} 
appear in the study of numerous physical phenomena where wave motion is
important, and hence methods for numerically calculating solutions to
these systems of partial differential equations have broad applications over
multiple disciplines. Complicating the development of numerical methods for
solving these systems is the fact that many solutions are not smooth, but
rather contain discontinuities that lead to computational difficulties. One
approach for accurately approximating such solutions is to use the
high-resolution finite volume method described in \cite{LeVeque1997,Leveque1}, which
focuses on calculating cell averages rather than a pointwise approximation
at grid points and uses Riemann solvers combined with limiters to avoid
non-physical oscillations.

Regardless of the accuracy of the particular algorithm being utilized, the
error associated with using a discrete grid to solve an originally continuous
problem is always present. This error can be reduced by increasing the
refinement of the mesh being utilized, but computation time considerations
limit the amount of refinement that is practical. Adaptive mesh refinement
(AMR) clusters grid points in areas of interest, such as discontinuities or
regions where the solution has a complicated structure. This allows 
refinement of grid cells in regions of interest, without spending
computational effort refining portions of the domain where the solution is smooth or
not of interest. A block-structured AMR algorithm developed to work in conjunction with
wave-propagation algorithms \cite{Berger1998} is available as the
AMRClaw package of Clawpack \cite{CLAWPACK}, and is also used in the GeoClaw
variant of Clawpack for modeling tsunamis, storm surge, and other geophysical flows (e.g.
\cite{BergerGeorgeLeVequeMandli:awr11,LeVequeGeorgeBerger2011,MandliDawson2014}).

A key component of any AMR algorithm is the criterion for deciding which grid cells
should be refined. In AMRClaw, one strategy is a Richardson 
error estimation procedure
that compares the solution on the existing grid with the solution on a coarser
grid, and refines cells where this error estimate is greater than a specified tolerance
\cite{Berger1998}.  Another approach is to simply flag cells where the
gradient of the solution (or an undivided difference of neighboring cell
values) is large.
For GeoClaw applications to tsunami and storm surge modeling, 
the approach often used is to flag cells
where the surface elevation of
the water is perturbed from sea level beyond some set tolerance.
These approaches would refine anywhere that the tolerance is exceeded,
irrespective of the fact that the area of interest may be only a subregion of
the full solution domain. 
To address this, recent versions of AMRClaw and GeoClaw also allow specifying
``refinement regions,'' space-time subsets of the computational domain where
refinement above a certain level can be either required or forbidden.  This
is essential in many GeoClaw applications where only a small region along the
coast (some community of interest) must be refined down to a very fine
resolution (often $1/3$ arcsecond, less than 10 meters) 
as part of an ocean-scale simulation.
These AMR regions can also be used to induce the code to follow only
the waves of interest, e.g. as a tsunami propagates across the ocean, but to
do so optimally often requires multiple attempts and careful examination of
how the solution is behaving, generally using coarser grid runs for guidance.
This manual guiding of AMR may also fail to capture some waves that are
important. For example, a portion of a tsunami wave may appear to be heading
away from the community of interest but later reflect off a distant
shoreline or underwater features, or edge waves may be excited that propagate
along the continental shelf for hours after the primary wave has passed.
These challenges in tsunami modeling were the original motivation for 
the work reported here and we include one example in \cref{sec:tsunami} to
illustrate our new approach.

For any problem where only a particular area of the total solution is of
interest, a method that allows specifically targeting and refining the
grid in regions that influence this area of interest would significantly
reduce the computational cost of finding a solution.  Other applications
where this could be very useful include earthquake simulation, for example, 
where again the desire might be to efficiently refine only the waves that
will reach a particular seismometer or community of interest.

In this paper we show how the time-dependent adjoint problem can be used in conjunction with
the AMR strategy already present in AMRClaw. We consider situations where the 
solution over a small subset of the full computational domain 
is of interest.  We first consider the case where we care about the
solution only at a single point in time, but then show that this
can be easily extended to the more typical case where we care about
a small spatial subset over a range of times, e.g. the full simulation time.
The code for all the examples presented is available online at \cite{adjointCode},
 and can be modified for use on other problems.

We only consider linear problems in this work, including examples in 1 and 2
space dimensions for variable-coefficient linear acoustics, and the tsunami
problem in \cref{sec:tsunami}.  For the latter problem we use GeoClaw, which
in general solves the nonlinear shallow water equations, but restrict our
attention to the above-mentioned application of tracking waves in the ocean
that will reach the region of interest.  Since a tsunami in the
ocean typically has a amplitude less than one meter, which is very small
compared to the ocean depth, these equations essentially reduce to the linear
shallow water equations and the adjoint equations linearized about the
ocean at rest is sufficient for our needs.  In \cref{sec:futurework} we make some
comments about extension to nonlinear problems.

Adjoint equations have been used computationally for many years in a variety 
of different fields, with wide ranging applications. A few examples 
include weather model tuning \cite{Hall1986}, aerodynamics 
design optimization \cite{GilesPierce2000,Jameson1988,KennedyMartins2013}, 
automobile aerodynamics \cite{Othmer2004}, and geodynamics
\cite{BungeHagelbergTravis2003}.  They have been used for
seismic inversion \cite{AkcelikBirosGhattas2002,TrompTapeLie2005}
and recently also applied to tsunami inversion \cite{Blaisea2013}.
The adjoint method 
has also been used for error estimation in the field of aerodynamics 
\cite{BeckerRannacher2001} and for general coupled time-dependent 
systems \cite{AsnerTavenerKay2012}. 
Various solution methods have been combined with
adjoint approaches, including Monte Carlo \cite{BuffoniCupini2001}, 
finite volume \cite{Mishra2013}, finite element 
\cite{AsnerTavenerKay2012}, and spectral-element 
\cite{TrompTapeLie2005} methods. 

The adjoint method has also been used to guide adaptive mesh 
refinement, typically by estimating the error in the calculation 
and using that to determine how to adjust the grid, e.g.
\cite{Carey2010,Nemec,NemecAftosmisWintzer2008,VendittiDarmofal2000}.
However, in these works steady state problems are solved and the
adjoint problem takes the form of another steady state problem,
typically yielding a large sparse system of linear equations to be
solved.  By contrast, in our time-dependent problem the adjoint
equation will also be time-dependent, and must be solved backwards
in time.

Computationally there are two different varieties of adjoint methods:
the discrete adjoint, where the PDE is first discretized and the adjoint 
system is derived by algorithmic differentiation of this discretization, and 
the continuous approach, where the adjoint of the original mathematical equation
is first derived and then discretized.
For some problems (e.g. to obtain precise error estimates or for optimizaiton
or control problems), it is crucial that the adjoint equations solved be the
adjoint of the discretized forward problem.
Li and Petzold \cite{LiPetzold2004} introduce a hybrid approach
combining aspects of discrete and continuous adjoint approaches for
computing sensitivities in time-dependent problems where AMR is
used, but do not discuss using the adjoint to guide AMR.

In our case we are primarily interested in identifying regions in space where
the grid for the forward problem should be refined and for this the
continuous adjoint approach appears to be sufficient and is much easier to
implement.  In \cref{sec:adjoint} we derive the adjoint equation for the PDE
of interest and show that it is also a hyperbolic PDE, which can be solved by
the same finite volume methods using Clawpack as are used for the forward
problem.  Determining the boundary conditions for the
adjoint PDE can be difficult for some problems and limits the ability to use
the continous adjoint approach (e.g. \cite{LiPetzold2004}), but for the
hyperbolic equations we consider we will see that these can be easily
obtained and implemented using the ghost cell methodology of Clawpack.

The discretization we use for the adjoint equation is not the adjoint
of the forward problem discretization and is not even solved on the same
computational grid.  In the work presented here, we 
only solve the adjoint problem once on a fixed grid and interpolate as needed
to the AMR grids being used for the forward problem (see \cref{sec:amradj}).  
More complicated approaches in which AMR is also applied to the adjoint problem 
are currently being studied.

\section{The Adjoint Equation}\label{sec:adjoint}

Suppose $q(x,t)$ is the solution to the time-depen\-dent linear equation 
(with spatially varying coefficients)
\begin{equation} \label{eq:qtAqx}
q_{t}(x,t) + A(x)q_{x}(x,t) = 0, \quad a\leq x \leq b, \quad t_0\leq t \leq
t_f
\end{equation} 
subject to some known initial conditions, 
$q(x,t_0)$, and some boundary conditions at $x=a$ and $x=b$. 
Here $q(x,t) \in \reals^m$ for a system of $m$ equations and
we assume $A(x) \in \reals^{m\times m}$ 
is diagonalizable with real eigenvalues at each $x$, so that
\cref{eq:qtAqx} is a hyperbolic system of equations.

Now suppose
we are interested in calculating the value of a functional
\begin{equation}\label{eq:J_general}
J = \int_a^b \varphi^T (x) q(x,t_f) dx
\end{equation} 
for some given $\varphi (x)$.
For example, if $\varphi (x) = \delta (x - x_0)$ then $J = q(x_0,t_f)$ is the
solution value at the point $x = x_0$ at the final time $t_f$. 
This is the situation we consider in this paper, with the delta function
smeared out around the region of interest for the computational approach.

If $\hat{q}(x,t)\in\reals^m$ is any other function then multiplying this by
\cref{eq:qtAqx} and integrating yields
\begin{equation}\label{eq:q_integral}
\int_a^b \int_{t_0}^{t_f}
\hat{q}^T(x,t)\left(q_{t}(x,t)+A(x)q_{x}(x,t)\right)dx\,dt = 0
\end{equation}
for any time $t_0 < t_f$,
then integrating by parts twice yields the equation 
\begin{equation}\label{eq:intbyparts}
\int_a^b   \left.\hat{q}^Tq\right|^{t_f}_{t_0}dx 
+ \int_{t_0}^{t_f} \left.\hat{q}^TAq\right|^{b}_{a}dt 
- \int_{t_0}^{t_f} \int_a^b  q^T\left(\hat{q}_{t} +
\left(A^T\hat{q}\right)_{x}\right)dx\,dt = 0.
\end{equation} 
By defining the adjoint equation, 
\begin{equation}\label{adjoint1}
\hat{q}_{t}(x,t) + (A^T(x)\hat{q}(x,t))_{x} = 0,
\end{equation}
setting $\hat{q}(x,t_f) = \varphi (x)$,
and selecting the appropriate boundary conditions for $\hat{q}(x,t)$ 
such that the integral in time vanishes (see below), we can eliminate all
terms from \cref{eq:intbyparts} except the first term, to obtain 
\begin{equation}\label{eq:q_equality}
\int_a^b \hat{q}^T(x,t_f)q(x,t_f) dx = \int_a^b
\hat{q}^T(x,t_0)q(x,t_0)dx.
\end{equation} 
Therefore, the integral of the inner product
between $\hat{q}$ and $q$ at the final time is equal 
to the integral at the initial time $t_0$:
\begin{equation}
J = \int_a^b\hat{q}^T(x,t_0)q(x,t_0)dx.\label{eq:J}
\end{equation}
Note that we can replace $t_0$ in \cref{eq:q_integral} with any 
$t$ so long as $t_0 \leq t \leq t_f$, which would yield \cref{eq:J} 
with $t_0$ replaced by $t$. 
 From this we observe that the locations where the inner product
$\hat{q}(x,t)^Tq(x,t)$ is large, for any $t$ with
$t_0 \leq t \leq t_f$, are the areas that will have a
significant effect on the inner product $J$.
These are the areas where the solution should be refined at time $t$.

To make use of this, we must first solve the adjoint equation
\cref{adjoint1} for $\hat q(x,t)$.
Note, however, that this requires using 
``initial'' data $\hat{q}(x,t_f)$, so
the adjoint problem must be solved backward in
time. The strategy used for this is discussed in \cref{sec:amradj}.

Note that it is also possible to use adjoint equations
to compute sensitivities of $J$ to changes in the input data, by using the
equation 
\begin{equation*}
\delta J = \int \hat{q}(x,t_0)\,\delta q(x,t_0)\,dx.
\end{equation*} 
This has led to the adjoint equations being utilized for system control 
in a wide variety of applications such as shallow-water wave control 
\cite{SandersKatopodes2000} and 
optimal control of free boundary problems \cite{Marburger2012}.
As a simple example, if we take
$\varphi(x) = \hat{q}(x,t_f)$ to be a delta function then this adjoint method 
approach would
provide us with the sensitivity of a single solution point to changes in the
 data $q(x,t_0)$. 
This is also useful in solving inverse problems and potential applications of
this approach in tsunami modeling are being studied separately.

\section{Numerical Methods}\label{sec:numerical}

We use Clawpack \cite{CLAWPACK} to solve both the forward problem and the
adjoint equation.
Clawpack utilizes logically rectangular grids, and each grid cell is viewed as
a volume over which cell averages of the solution variables are calculated. At
each time step, the interface between two cells is treated as a Riemann
problem that is solved to compute the waves propagating into the grid
cells. The {\em Riemann problem} is simply the initial value problem together with piecewise
constant data, which is determined by the cell averages of the dependent
variables. The cell averages are then updated by the waves propagating into
the grid cell from each cell edge. These Riemann problems are solved using
{\em Riemann solvers,} or more commonly {\em approximate Riemann solvers,} which
vary depending on the specific equations being solved. In two-dimensional
problems a {\em transverse Riemann problem} must also be solved, in which the
waves moving normal to a cell edge are split in the transverse direction and
are used to modify the cell averages in adjacent rows of grid cell. This
approach is presented in more detail in \cite{LeVeque1997}, and various
examples can be found in Chapters 20-22 of \cite{Leveque1}.

Note that the forward problem \cref{eq:qtAqx} is not in conservation
form (unless $A$ is constant), and the examples we consider below
have this form.  The adjoint equation \cref{adjoint1}, however, is
a hyperbolic problem in conservation form.  Conversely, if we started
with a linear problem in conservation form then the adjoint problem
obtained via integration by parts would be in non-conservation form.
Either way, it is important to use methods that can handle
non-conservative hyperbolic systems, a feature of the general
wave-propagation methods used in Clawpack.
Another approach available in Clawpack is the f-wave formulation for
conservative equations with spatially varying flux functions, required for
the nonlinear shallow water equations over varying topography and originally 
presented for shallow water equations in
\cite{BaleLeVequeMitranRossmanith2002}. In the f-wave approach, it is the
difference in the flux normal to the interface between two cells that is split
into waves, and these waves are then used to update the adjacent cells.
GeoClaw uses this f-wave formulation to solve an augmented Riemann problem
that also robustly handles wetting and drying in inundation zones
\cite{dgeorge:jcp}.
Using this f-wave approach to solve the adjoint equations
for shallow water problems 
has been shown to yield accurate results \cite{SandersBradford2002}.
With the aim of incorporating the adjoint method
code developed for AMR into GeoClaw, the examples provided in this paper also
utilize f-wave formulations for the adjoint even when the 
problem does not pertain to
geophysical flows.

The mesh refinement used in AMRClaw, which is incorporated into the Clawpack
software, consists of multiple levels of nested patches,
each of which is generally refined in both time and space to preserve the
stability of the finite volume method. This means that for each time step
taken on a coarser grid, multiple steps must be taken on finer grids. At each
time step it is necessary to fill in {\em ghost cell} values around each finer
grid in order to provide accurate boundary conditions for the time step. For
each ghost cell, the values come from either a neighboring grid at the same
level, if such a grid exists, or by interpolating in both time and space from
the values in the underlying coarse grid. This requires that the code be
organized to advance the time step on the coarsest grid first, and then update
the necessary number of time steps on finer grids. This same procedure is used
recursively at all the grid levels, meaning that the finest patches are
updated last.

Grids are refined by flagging cells were the resolution is determined to be
insufficient, and then clustering the flagged cells into rectangular
refinement patches. An integral part of the refinement algorithm, which is
presented in \cite{BergerOliger1984} and \cite{BergerColella1989}, is
clustering the flagged cells into grid patches that limit the number of
unflagged cells contained in the grids while also not introducing too many
separate patches using an algorithm of Berger and Rigotsis 
\cite{BergerRigoutsos1991}. The grids are also constrained to follow a nesting
criterion, so any grid patch at level $n$ must be contained in a grid patch at
level $n -1$. Every few time steps the features in the solution that require
refinement will have moved, which requires the grid patches to be recalculated
based on newly flagged cells. Based on this, the methodology utilized for
selecting which cells to flag will have a significant impact on both the
accuracy of the results and the time required for the computation to run.

\section{Combining Mesh Refinement and the Adjoint Problem}\label{sec:amradj}
By developing a strategy for taking advantage of the sensitivity information
provided by adjoint methods and incorporating it into an AMR algorithm, it is
possible to significantly reduce the computational time required to find the
solution. As examples, we present acoustics problems in both one and
two dimensions, and a tsunami modeling case.

In utilizing the adjoint solution to guide mesh refinement for the original
problem, which we will refer to as the ``forward problem,'' a fairly coarse
grid solution to the adjoint equation 
has produced excellent results. Exploring the implications and
benefits of allowing for mesh refinement when solving 
the adjoint problem is a possible
direction for future work. For the examples shown later in the paper, the
adjoint problem is solved on a grid that is not refined. 

Consider the one dimensional problem \cref{eq:qtAqx} and recall that
the adjoint equation \cref{adjoint1} has the form
\begin{align*}
\hat{q}_t + \left( A^T(x) \hat{q}\right)_x  = 0,
\end{align*}
where the initial condition for $\hat{q}$ is given at the final time,
$\hat{q}(x,t_f) = \varphi (x)$, and is selected to highlight the
impact of the solution on some region of interest. Clawpack is designed to
solve equations forward in time, so we consider the function
\begin{align*}
\tilde{q}(x,t) &= \hat{q}(x,t_f - t).
\end{align*}
This gives us the new problem
\begin{align*}
&\tilde{q}_t - \left( A^T(x) \tilde{q}\right)_x = 0 
&&x \in [a,b], \hspace{0.1in}t > 0\\
&\tilde{q}(a,t) = \hat{q}(a,t_f - t)  &&0 \leq t \leq t_f - t_0 \\
&\tilde{q}(b,t) = \hat{q}(b,t_f - t) &&0 \leq t \leq t_f - t_0
\end{align*}
with initial condition $\tilde{q}(x,0) = \varphi (x) $. This
problem is then solved using the Clawpack software. Snapshots of this solution
are saved at regular time intervals, $t_0, t_1, \cdots, t_N$. After the
solution is calculated, snapshots of the adjoint solution are retrieved by
simply setting
\begin{align*}
\hat{q}(x,t - t_n) = \tilde{q}(x,t_n)
\end{align*}
for $n = 0, 1, \cdots, N$. 

With the adjoint solution in hand, we now turn to the forward problem. As
refinement occurs in space for the forward problem, maintaining the stability
of the finite volume method requires that refinement must also occur in time.
Therefore, as the forward solution is refined, solution data for the adjoint
problem is no longer available at the corresponding times (since only
snapshots at regular time intervals of the adjoint solutions were saved) or
locations (since the adjoint solution was calculated on a coarse grid). To
address this issue, the solution for the adjoint problem at the necessary
locations is approximated using bilinear interpolation from the data present
on the coarser grid. To be conservative, when considering the forward problem
at time $t$ with
\begin{align*}
t_n \leq t \leq t_{n+1},
\end{align*}
 both $\hat{q}(x,t_n)$ and $\hat{q}(x,t_{n+1})$ are taken
into account.

When solving the forward problem we are then able to take the inner product
between the current time step in the forward problem and the two corresponding
time steps in the adjoint problem in order to determine which areas in the
forward wave are going to impact the region of interest. These areas are then
flagged for refinement, and the next time step is taken.

\section{Numerical Results}
We present several examples using linear acoustics and tsunami modeling.
The timing
benefits of using the adjoint method to guide mesh refinement are shown
through the two dimensional examples in \cref{sec:acou2dperf,sec:tsunamiperf}.

\subsection{Acoustics In One Space Dimension}
Clawpack contains AMR implementations in 2D and 3D but not for
one-dimensional problems. However, we start by describing our strategy in 1D
since the adjoint approach is best visualized through these examples. 
Consider the linear acoustics equations in one dimension in a piecewise
constant medium, with wall boundary conditions on both the left and the right:
\begin{align*}
p_t(x,t) + K(x) u_x(x,t) &= 0 &&x \in [a,b], \hspace{0.1in}t > t_0,\\
{\rho (x)}u_t(x,t) + {p_x(x,t)} &= 0 &&x \in [a,b],
\hspace{0.1in}t > t_0,\\
u(a,t) = 0, \hspace{0.1in} u(b,t) &= 0 &&t  \geq t_0.
\end{align*}
Here, $p$ is the pressure, $u$ is the velocity, $K$ is the bulk modulus, and
$\rho$ is the density. Setting 
\begin{align*}
A(x) = \left[ \begin{matrix}
0 & K(x) \\
1/\rho(x) & 0
\end{matrix}\right], \hspace{0.3in}
q(x,t) = \left[\begin{matrix}
p(x,t) \\ u(x,t)
\end{matrix}\right],
\end{align*}
gives us the equation $q_t(x,t) + A(x)q_x(x,t) = 0$. Suppose that 
we are interested in the accurate estimation of the pressure in the 
interval $1.8 < x < 2.3$. Setting $J = \int_{1.8}^{2.3}p(x,t_f)\,dx$,
the problem then requires that
\begin{align}
\varphi (x) = \left[ \begin{matrix}
I (x) \\ 0
\end{matrix}
\right], \label{eq:phi_1d}
\end{align}
where
\begin{align}
I (x) = \left\{
     \begin{array}{lr}
       1 & \hspace{0.3in}\textnormal{if } 1.8 < x < 2.3\\
       0 & \textnormal{otherwise.}\hspace{0.325in}
     \end{array}
   \right. \label{eq:delta_1d}
\end{align}
Define
\begin{align*}
\hat{q}(x,t_f) = \left[ \begin{matrix}
\hat{p}(x,t_f) \\ \hat{u}(x,t_f)
\end{matrix}\right] = \varphi(x),
\end{align*}
and note that \cref{eq:q_integral} holds for this problem.
If we define the adjoint problem
\begin{align*}
\hat{q}_{t} + \left(A^T(x)\hat{q}\right)_{x} &= 0 &&x \in [a,b],
\hspace{0.1in}t > t_0\\
\hat{u}(a,t) = 0, \hspace{0.1in} \hat{u}(b,t) &= 0 &&t  \geq t_0,
\end{align*}
then we are left with \cref{eq:q_equality},
which is the expression that allows us to use the inner product of the
adjoint and forward problems at each time step to determine what regions will
influence the point of interest at the final time. 

As an example, consider a linear system with piecewise constant coefficients 
where $t_0 = 0$, $a = -5$, $b = 3$, $K = 1$ and 
\begin{align*}
\rho (x) = \left\{
     \begin{array}{lr}
       1 & \hspace{0.3in}\textnormal{if } x < 0\hspace{0.04in}\\
       4 & \textnormal{if } x > 0,
     \end{array}
   \right.
\end{align*}
giving a jump in sound speed from $c = 1$ on the left to $c = 0.5$ 
on the right. As initial
data for $q(x,t)$ we take a Gaussian hump in pressure in the left region, and
a zero velocity. The initial hump in pressure is given by
\begin{align*}
p(x,0) = \exp \left( - \beta(x + 2)^2\right)
\end{align*}
with $\beta = 50$. As time progresses, the hump splits into 
equal left-going and right-going waves which interact with the walls and 
the interface giving both reflected and transmitted waves. 

\begin{figure}[ht]
 \centering
  \subfigure[Results for $t = 0.0$ seconds.]{
  \includegraphics[width=0.475\textwidth]{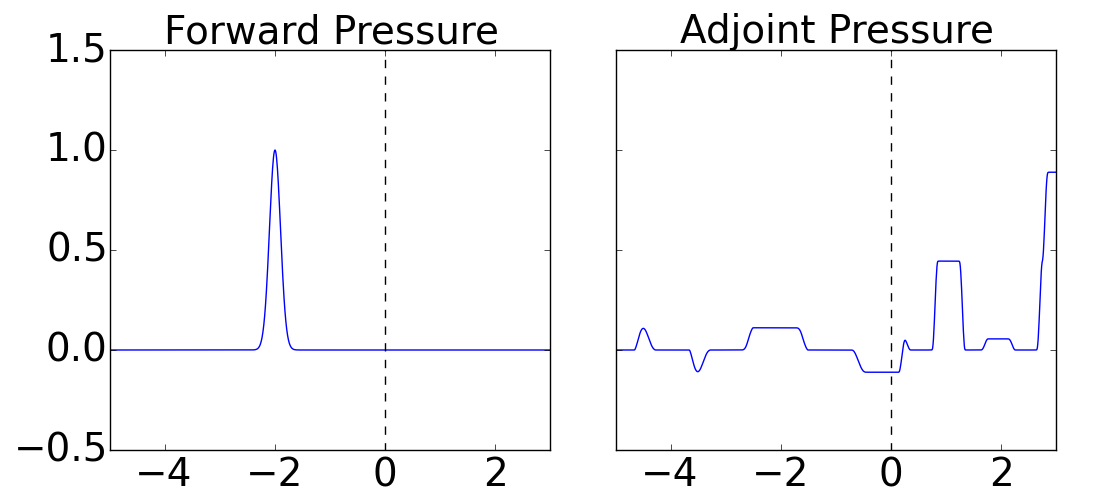}
   \label{fig:1dAcoustics_1}
   }
  \subfigure[Results for $t = 1.5$ seconds.]{
  \includegraphics[width=0.475\textwidth]{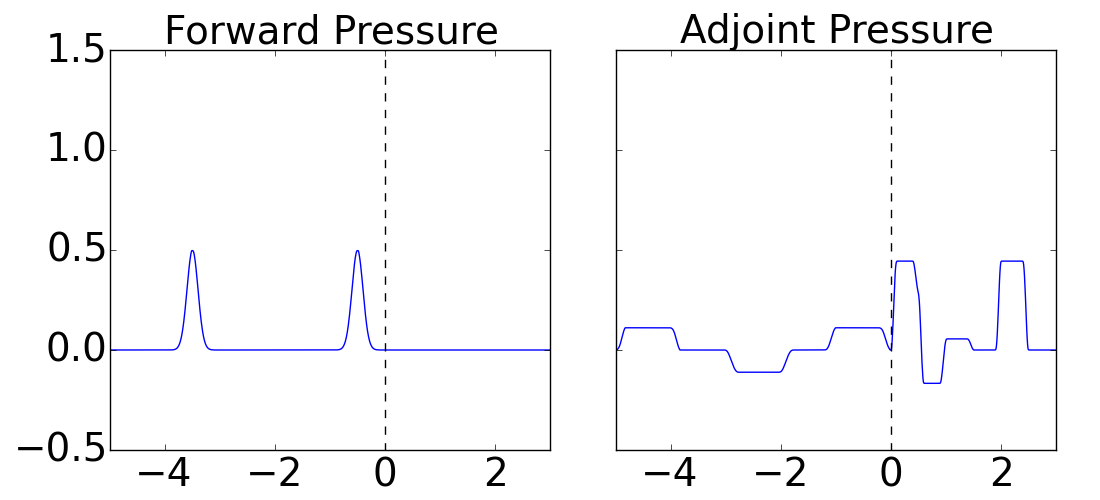}
   \label{fig:1dAcoustics_2}
   }
  \subfigure[Results for $t = 2.5$ seconds.]{
  \includegraphics[width=0.475\textwidth]{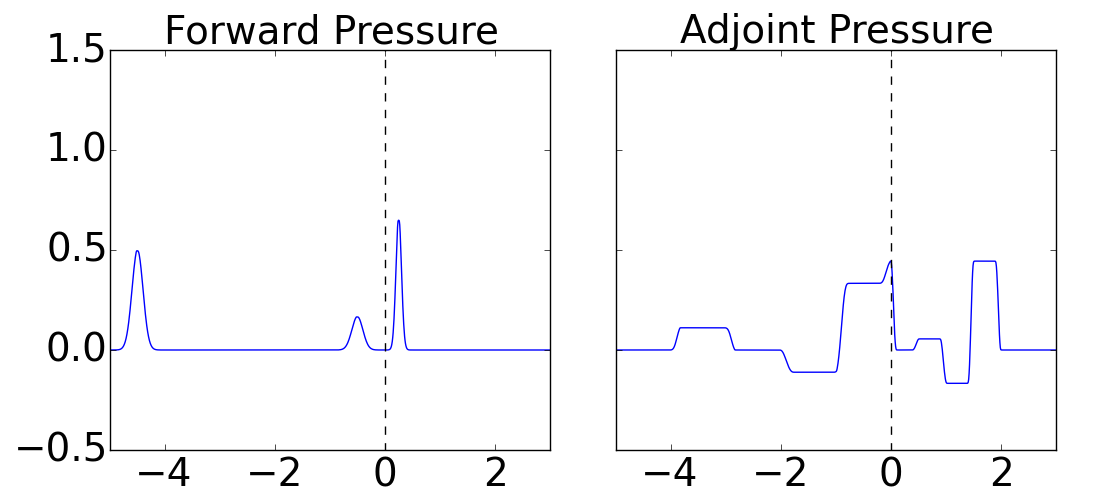}
   \label{fig:1dAcoustics_3}
   }
  \subfigure[Results for $t = 15.0$ seconds.]{
  \includegraphics[width=0.475\textwidth]{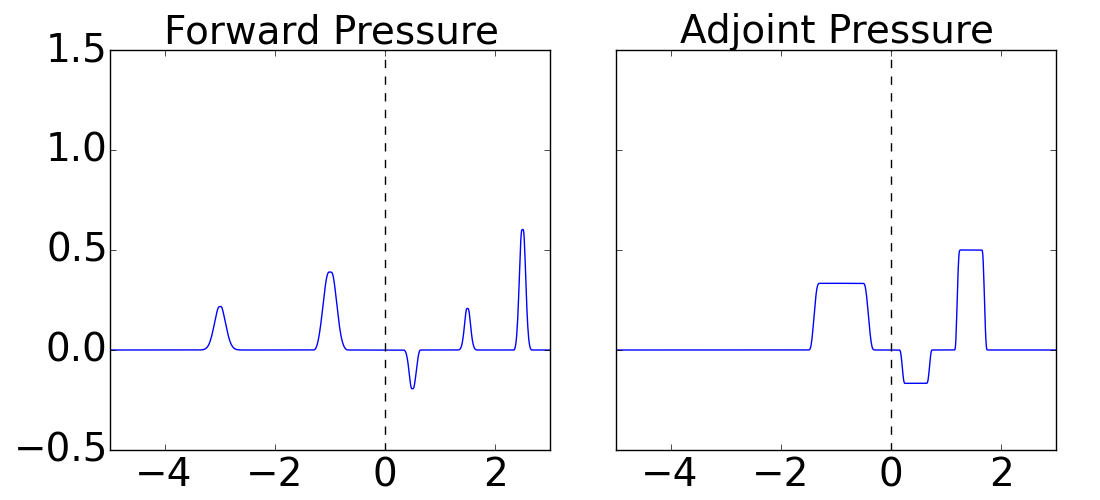}
   \label{fig:1dAcoustics_4}
   }
  \subfigure[Results for $t = 19.0$ seconds.]{
  \includegraphics[width=0.475\textwidth]{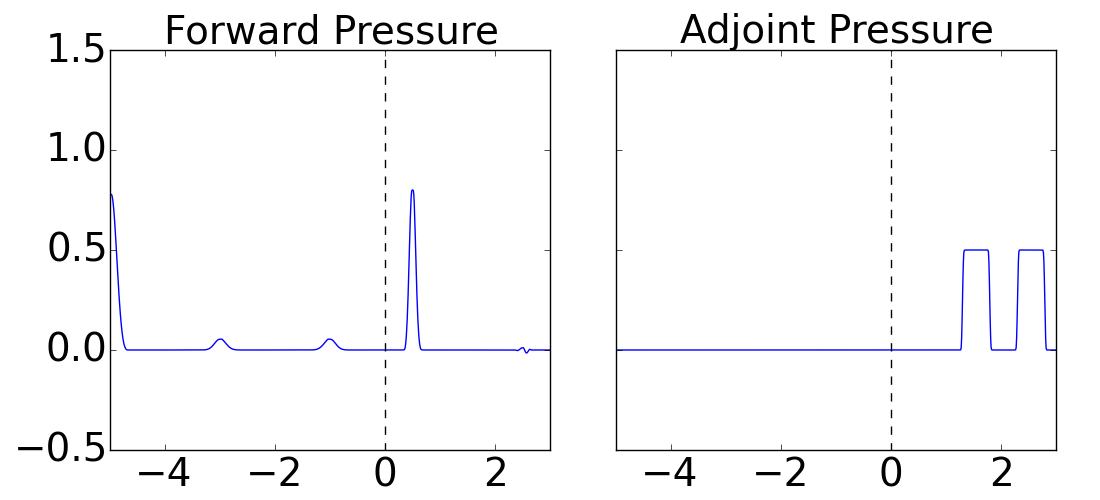}
   \label{fig:1dAcoustics_5}
   }
  \subfigure[Results for $t = 20.0$ seconds.]{
  \includegraphics[width=0.475\textwidth]{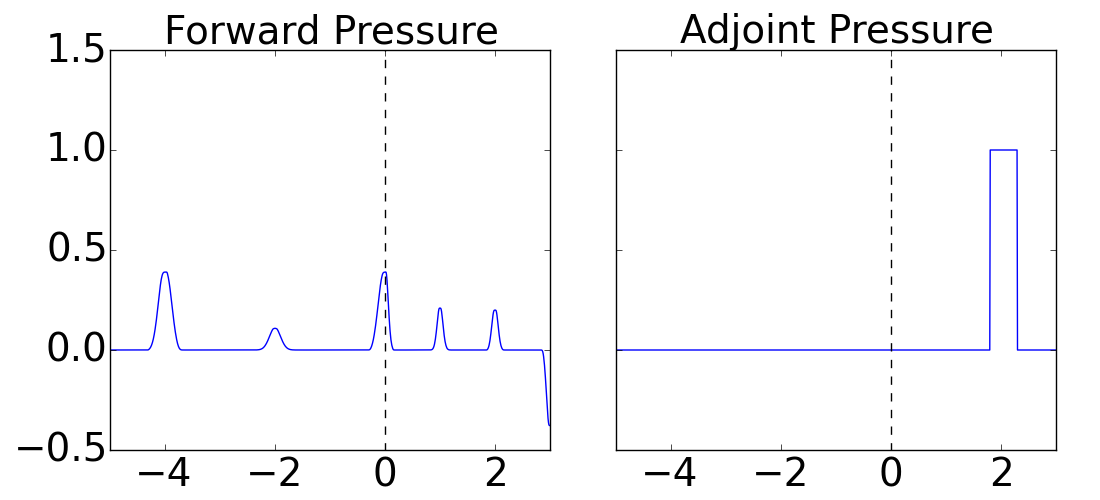}
   \label{fig:1dAcoustics_6}
   }
 \caption[Optional caption for list of figures]{%
  Computed results for one-dimensional acoustics problem.
  The $y$-axis is the same for both plots in each subfigure.
  The adjoint pressure evolves backward in time from (f) to (a).}
   \label{fig:1dAcoustics}
\end{figure}

As the ``initial'' data for $\hat{q}(x,t_f) = \varphi (x)$ 
we have a square pulse in pressure, which was described above in
Equations \cref{eq:phi_1d} and \cref{eq:delta_1d} at the final time.
As time progresses backwards, the pulse splits into equal left-going and right-going
waves which interact with the walls and the interface giving both reflected
and transmitted waves. \cref{fig:1dAcoustics} shows the Clawpack 
results for the pressure of both the forward and adjoint solutions at
six different times, on a grid with 1000 grid points and no mesh refinement. 
Both of these problems are run using $t_f = 20$, so the forward solution is
run from $t= 0$ to $t = 20$ and the adjoint solution is run from $t= 20$ to $t
= 0$.

\begin{figure}[ht]
 \centering
  \subfigure[The 1-norm
of $q(x,t)$ and $\hat{q}(x,t)$ are shown in the left and right figures
respectively, showing all waves.]{
  \includegraphics[width=0.475\textwidth]{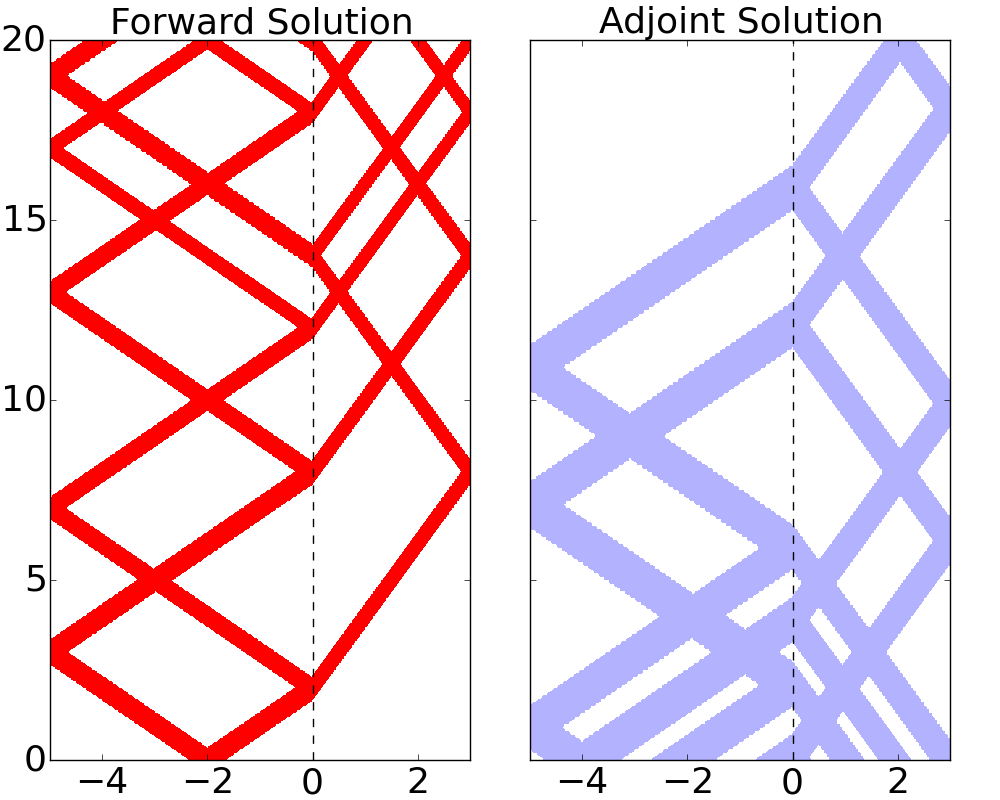}
   \label{fig:1dAcoustics_norms}
   }
 \subfigure[The 1-norm 
of $q(x,t)$ and $\hat{q}(x,t)$ are overlayed in the left figure. Their inner
product is shown in the right, indicating where the forward problem should be
refined when only the solution $q(x,t)$
for $1.8<x<2.3$ at time $t=20$ is of interest.]{
  \includegraphics[width=0.475\textwidth]{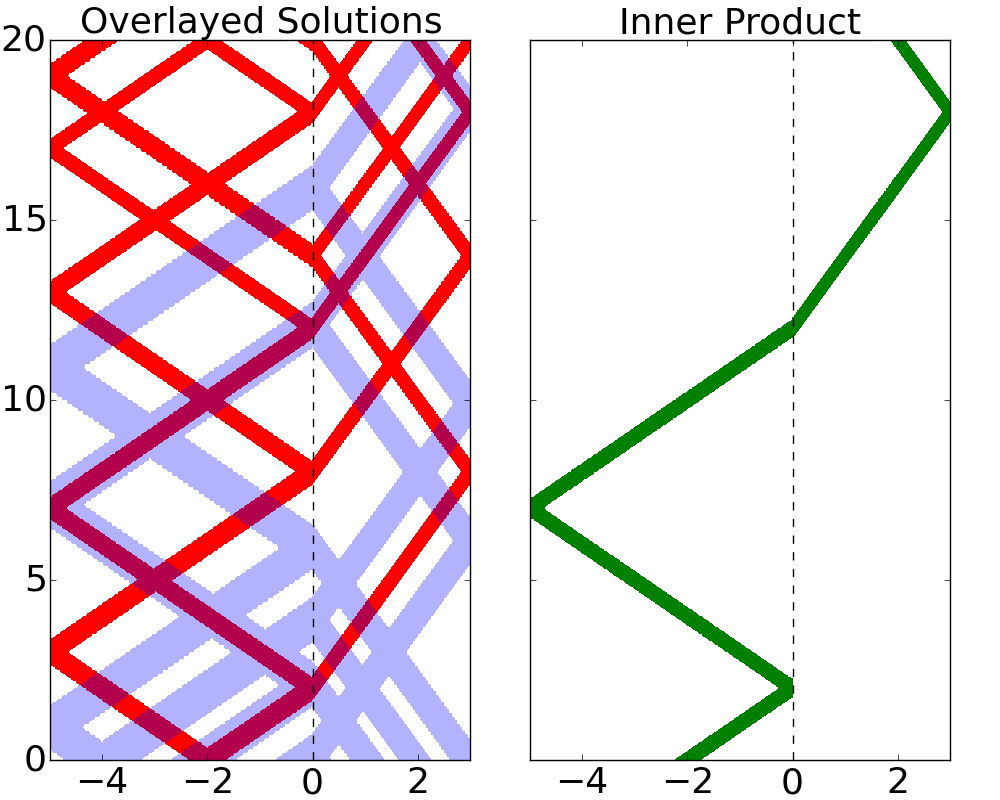}
   \label{fig:1dAcoustics_point}
   }
 \label{fig:1dAcoustics_xtplane}
 \caption[Optional caption for list of figures]{%
  Computed results for one-dimensional acoustics from \cref{fig:1dAcoustics}, 
  shown in the $x$-$t$ plane.
  The $y$-axis is the same for both plots in each subfigure.}
\end{figure}

To better visualize how the waves are moving through the domain, it is helpful
to look at the data in the $x$-$t$ plane as shown in \cref{fig:1dAcoustics_norms}. 
For \cref{fig:1dAcoustics_norms}, the
horizontal axis is the position, $x$, and the vertical axis is time. The left
plot shows in red the locations where the $1-$norm of $q(x,t)$ is greater than or
equal to $0.1$. The right plot shows in blue the locations where the $1-$norm
of $\hat{q}(x,t)$ is greater than or equal to $0.1$.

\subsubsection{Single Point In Time}
Suppose that we are interested in the solution $q(x,t)$ in the region given
by $1.8 < x < 2.3$ at the time $t = 20$. We have already computed the adjoint
solution, and have shown that by taking the inner product between the adjoint
solution and the forward solution it is possible to identify the regions that
will actually influence our region of interest at the final time. In the left
side of \cref{fig:1dAcoustics_point} we have overlayed the forward and 
adjoint solutions, and 
viewing the data in the $x$-$t$ plane makes it fairly clear which parts 
of the wave from the forward
solution actually effect our region of interest at the final time. 

\Cref{fig:1dAcoustics_point}
also shows, on the right, the locations where the inner product between the
forward and adjoint solution is greater than or equal to 0.1 as time
progresses. At each time step this is clearly identifying the regions in the
computational domain that we are interested in. If we were using adaptive mesh
refinement, these areas identified by where $\hat{q}^T(x,t)q(x,t)$ 
exceeds some tolerance are precisely the areas we would flag for 
refinement. Note that a mesh refinement strategy 
based on wherever the $1-$norm of $q(x,t)$ is large would result in refinement of
many areas in the computational domain that will have no effect on our area of
interest at the final time (all the red regions in the left plot of 
\cref{fig:1dAcoustics_norms}).

\subsubsection{Time Range}
Suppose that we are interested in the accurate estimation of the pressure in
the interval given
by $1.8 < x < 2.3$ for the time range $t_s \leq t \leq t_f$, where 
$t_s = 18$ and $t_f = 20$. Define $\hat q(x,t; \overline{t})$ as the adjoint based
on data $\hat q(x,\overline{t}) = \varphi(x)$.  Then for each $\hat t$ in the interval
$[t_s,t_f]$, we need to consider the inner product of $q(x,t)$ with $\hat q(x,t;
\hat t)$. Note that since the adjoint is autonomous 
in time, $\hat q(x,t; \hat t) = \hat q(x, t_f-\hat t+t; t_f)$. Therefore, we must 
consider the inner product 
\begin{equation*}
\hat q^T(x, t_f-\hat t+t; t_f)q(x,t)
\end{equation*}
for $\hat t \in [t_s,t_f]$. Since we are in fact only concerned when this 
inner product is greater than some tolerance, we can simply consider
\begin{equation*}
\max\limits_{t_s \leq \hat t \leq t_f} \hat q^T(x, t_f-\hat t+t; t_f)q(x,t)
\end{equation*}
and refine when this maximum inner product is above the given tolerance. 
Define $\tau = t_f-\hat t+t$. Then this maximum inner product can be 
rewritten as
\begin{align*}
\max\limits_{T\leq \tau \leq t} \hat{q}^T(x,\tau;t_f)q(x,t)
\end{align*}
where $T = \min(t+t_f-t_s,~ t_f)$.

\begin{figure}[h!]
  \centering
    \includegraphics[width=0.475\textwidth]{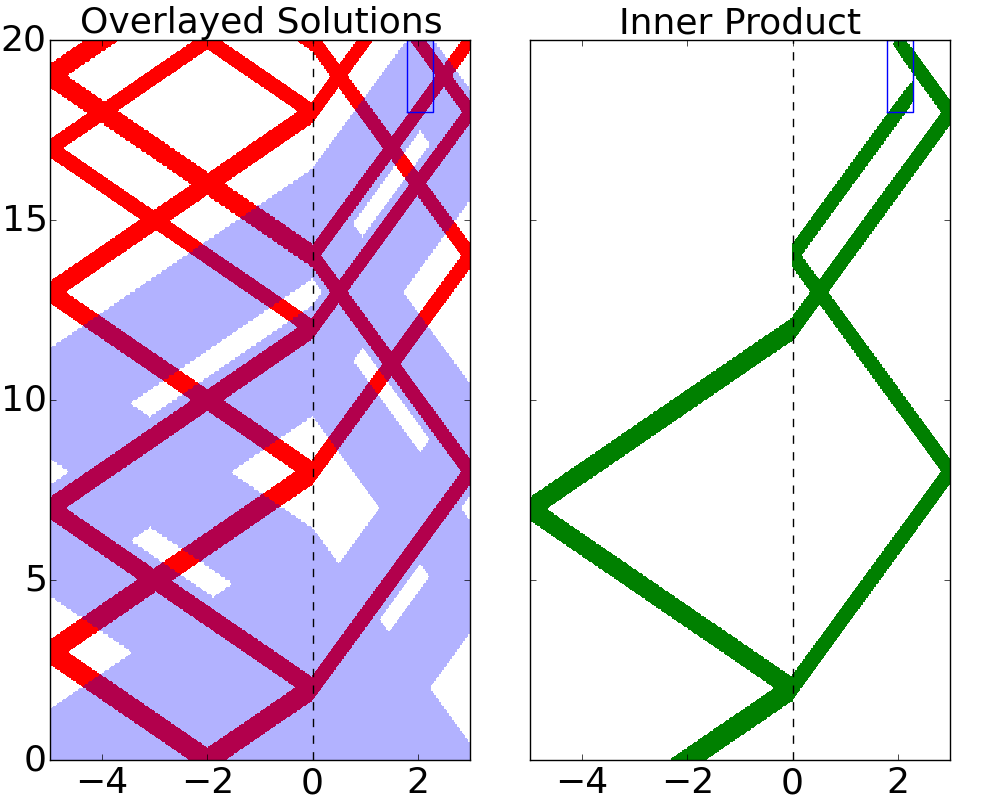}
    \caption{Computed results for one-dimensional acoustics in the case
when the solution $q(x,t)$ for $1.8<x<2.3$ and $18\leq t \leq 20$ is of
interest, as indicated by the box in each $x$--$t$ plot.
Left: 1-norm for 
    $q(x,t)$ and $\hat{q}(x,\tau)$ for $18\leq \tau \leq 20$. Right: Max inner product over the given 
    time range, showing the additional wave that must be refined.
    The time axis is the same for both plots.}
    \label{fig:1dAcoustics_timerange}
\end{figure}

We now drop the cumbersum notation $\hat q(x,t; t_f)$ in favor of the simpler 
$\hat{q}(x,t)$ with the understanding that the adjoint is based on the data
$\hat q(x,t_f) = \varphi(x)$. The left side of 
\cref{fig:1dAcoustics_timerange} shows in red the locations where the 1-norm 
of $q(x,t)$ is greater than or equal to 0.1 and in blue the locations where 
the 1-norm of $\hat{q}(x,\tau)$ is greater than or equal to 0.1, 
for $T \leq \tau \leq t$. \Cref{fig:1dAcoustics_timerange} also shows, 
on the right, the locations where the maximum inner product equation 
given above exceeds the tolerance of 0.1. 
At each time step this is clearly identifying the regions in the computational 
domain that will influence the interval of interest in the given time range. 
Therefore, the regions identified by this method are the regions that should
be refined at each time step.

\subsection{Acoustics in Two Space Dimensions}\label{sec:acoustics2d}
In two dimensions the linear acoustics equations are 
\begin{align*}
p_t(x,y,t) + K(x,y) \left(u_x(x,y,t) + v(x,y,t)_y\right) &= 0 &&x
\in [a,b],  y \in [\alpha,\beta], t > t_0\\
\rho(x,y) u_t(x,y,t) + {p_x(x,y,t)} &= 0 &&x \in [a,b],
 y \in [\alpha,\beta], t > t_0 \\
\rho(x,y) v_t(x,y,t) + {p_y(x,y,t)} &= 0 &&x \in [a,b],
 y \in [\alpha,\beta], t > t_0.
\end{align*} 
Setting 
\begin{align*}
A(x,y) = \left[ \begin{matrix}
0 & K & 0 \\
1/\rho& 0 & 0 \\
0 & 0 & 0
\end{matrix}\right],\hspace{0.1in}
B(x,y) = \left[ \begin{matrix}
0 & 0 & K \\
0 & 0 & 0 \\
1/\rho& 0 & 0
\end{matrix}\right], \hspace{0.1in}
q(x,y,t) = \left[\begin{matrix}
p \\ u \\ v
\end{matrix}\right],
\end{align*}
gives us the equation $q_t(x,y,t) + A(x,y)q_x(x,y,t) + B(x,y)q_y(x,y,t) = 0$. 

Suppose that we are interested in the accurate estimation of the 
pressure in the area defined by a square centered about 
$(x,y) = (3.56,0.56)$.
Setting 
\begin{align*}
J = 2\int_{3.32}^{3.8}\int_{0.32}^{0.8}p(x,y,t_f)dy\,dx,
\end{align*}
the problem then requires that
\begin{align}
\varphi (x,y) = \left[ \begin{matrix}
2 I(x,y) \\ 0 \\ 0
\end{matrix}
\right], \label{eq:phi_2d}
\end{align}
where
\begin{align}
I (x,y) = \left\{
     \begin{array}{lr}
       1 & \hspace{0.3in}\textnormal{if } 3.32 \leq x \leq 3.8
        \textnormal{ and } 0.32 \leq y \leq 0.8,\\
       0 & \textnormal{otherwise.}\hspace{1.61in}
     \end{array}
   \right. \label{eq:delta_2d}
\end{align}
Note that the coefficient of 2 in the definition of the functional $J$ is chosen for 
convenience. This coefficient can be chosen to be any value, provided that the 
refinement tolerance later on is adjusted accordingly.
Define
\begin{align*}
\hat{q}(x,y,t_f) = \left[ \begin{matrix}
\hat{p}(x,y,t_f) \\ \hat{u}(x,y,t_f) \\ \hat{v}(x,y,t_f)
\end{matrix}\right] = \varphi(x,y),
\end{align*}
and note that for any time $t_0 < t_f$ we have
\begin{align*}
\int_{t_0}^{t_f}\int_a^b\int_{\alpha}^{\beta}  \hat{q}^T\left( q_t + A(x,y)q_x +
B(x,y)q_y\right) dy\,dx\,dt = 0.
\end{align*}
Integrating by parts yields the equation
\begin{align}
\left. \int_a^b  \int_{\alpha}^{\beta}  \hat{q}^Tq\,dy\,dx \right|^{t_f}_{t_0}
&+ \left. \int_{t_0}^{t_f}\int_{\alpha}^{\beta}  \hat{q}^TA(x,y)q\,dy\,dt \right|^{b}_{a}
+ \left. \int_{t_0}^{t_f}\int_{a}^{b}  
\hat{q}^TB(x,y)q\,dx\,dt \right|^{\beta}_{\alpha}\nonumber \\
- &\int_{t_0}^{t_f}\int_a^b\int_{\alpha}^{\beta}  q^T\left(\hat{q}_{t} +
\left(A^T(x,y)\hat{q}\right)_{x} + \left(B^T(x,y)\hat{q}\right)_{y}\right)
dy\,dx\,dt = 0.\label{eq:acousticseqn_2d}
\end{align}
Note that if we can define an adjoint problem such that all but the first term 
in this equation vanishes then we are left with
\begin{align*}
\int_a^b \int_{\alpha}^{\beta}  \hat{q}^T(x, y, t_f)q(x, y, t_f) dy\,dx
= \int_a^b \int_{\alpha}^{\beta} \hat{q}^T(x, y, t_0)q(x, y, t_0)dy\,dx,
\end{align*}
which is the expression that allows us to use the inner product of the
adjoint and forward problems at each time step to determine what regions will
influence the point of interest at the final time. 
The adjoint problem required to accomplish this depends on the boundary conditions 
of the forward problem in question. Two different boundary conditions for the 
forward problem and the corresponding adjoint problems are discussed 
further below.

As an example,
consider $t_0 = 0$, $a = -4$, $b = 8$, $\alpha = -1$, $\beta = 11$, $K = 4$ 
and $\rho = 1$. As initial
data for $q(x,y,t)$ we take a smooth radially symmetric hump in pressure, and
a zero velocity in both $x$ and $y$ directions. The initial hump in pressure is given by
\begin{align}
p(x,y,0) = \left\{
     \begin{array}{lr}
       3 + \cos \left((\pi \left(r - 0.5\right)/w\right) 
       & \hspace{0.3in}\textnormal{if } \left| r - 0.3 \right| \leq w,\\
       0 & \textnormal{otherwise.}\hspace{0.353in}
     \end{array}
   \right. \label{eq:2d_acoustics_p}
\end{align}
with $w = 0.15$ and $r = \sqrt{\left(x-0.5\right)^2 +\left( y - 1\right)^2}$. 
As time progresses, this hump in pressure will radiate outward symmetrically.

Since the time interval of interest is given by $t_s \leq t \leq t_f$, 
at time $t$ we will refine the computational mesh where 
\begin{align*}
\max\limits_{T \leq \tau \leq t} \hat{q}^T(x,y,\tau ) q(x,y,t)
\end{align*}
is above some tolerance, with $T = \min (t + t_f - t_s, 0)$. For the examples
below this tolerance is set to $0.02$. 

Note that adaptive mesh refinement was used for the forward problem in these 
three examples, and the black outlines in each figure are the edges of different 
refinement levels. In each example, three levels of refinement are allowed
for the forward problem,
with a refinement ratio of 2 from each grid to the next. 
The adjoint problems were solved
on a $50 \times 50$ grid with no mesh refinement,
which is also the resolution of the coarsest level for the 
forward problems.

\subsubsection{Wall boundary conditions for a time point of interest}\label{sec:acoustics_timepoint}
Suppose that $t_f = 1.5$ seconds and that we are interested in the accurate calculation
of the pressure only at the final time (so $t_s = t_f = 1.5$). If we have the boundary conditions 
\begin{align}
&u(a,y,t) = 0, \hspace{0.1in} u(b,y,t) = 0 &&t  \geq 0, \nonumber \\
&v(x,\alpha,t) = 0, \hspace{0.1in} v(x,\beta,t) = 0 &&t
  \geq 0, \label{eq:acoustics2d_wallsBC}
\end{align}
then all but the first term in Equation \cref{eq:acousticseqn_2d} vanish if we 
define the adjoint problem
\begin{align}
&\hat{q}_{t} + \left(A^T(x,y)\hat{q}\right)_{x} + \left(B^T(x,y)
\hat{q}\right)_{y} = 0
&&x \in [a,b], y \in [\alpha, \beta], t > 0 \nonumber \\
&\hat{u}(a,y,t) = 0, \hspace{0.1in} \hat{u}(b,y,t) = 0 &&t  \geq 0  
\label{eq:acoustics2d_adjointwalls} \\
&\hat{v}(x,\alpha,t) = 0, \hspace{0.1in} \hat{v}(x,\beta,t) = 0
&&t  \geq 0. \nonumber 
\end{align}
\begin{figure}[ht]
 \centering
  \subfigure[Results at $t = 0$.]{
  \includegraphics[width=0.47\textwidth]{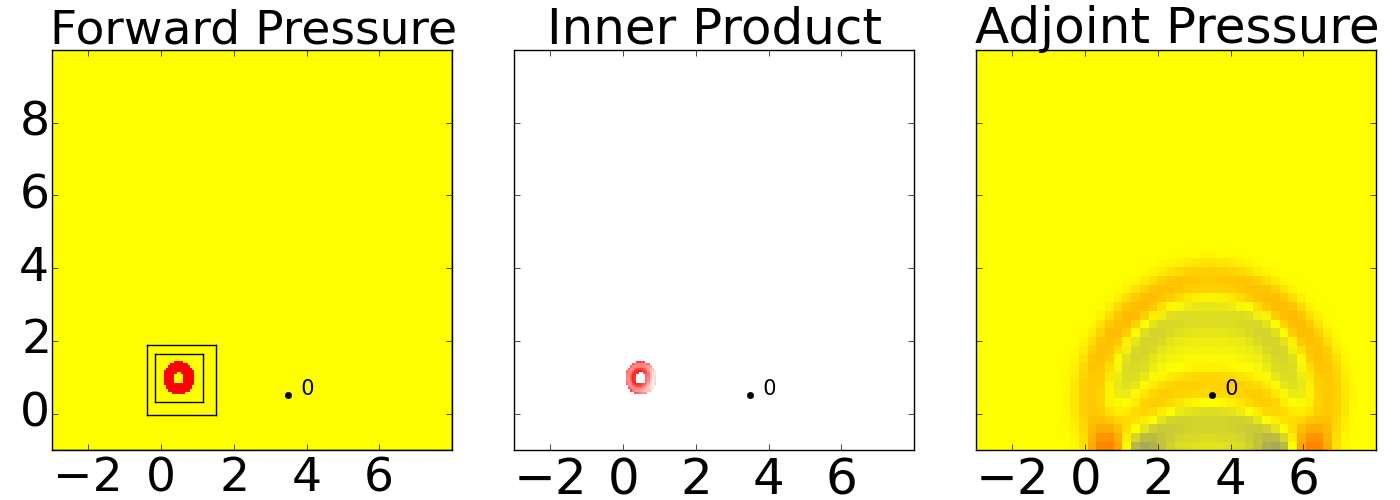}
   \label{fig:2dAcoustics_tp1}
   }
   \subfigure[Results at $t = 0.5$.]{
  \includegraphics[width=0.47\textwidth]{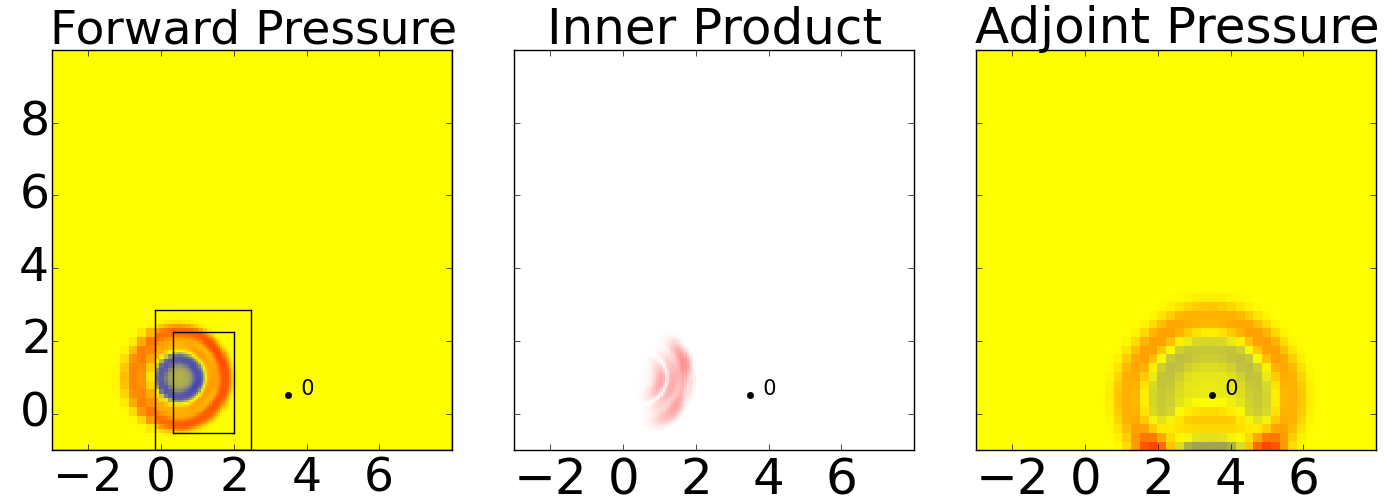}
   \label{fig:2dAcoustics_tp2}
   }
  \subfigure[Results at $t = 1.0$.]{
  \includegraphics[width=0.47\textwidth]{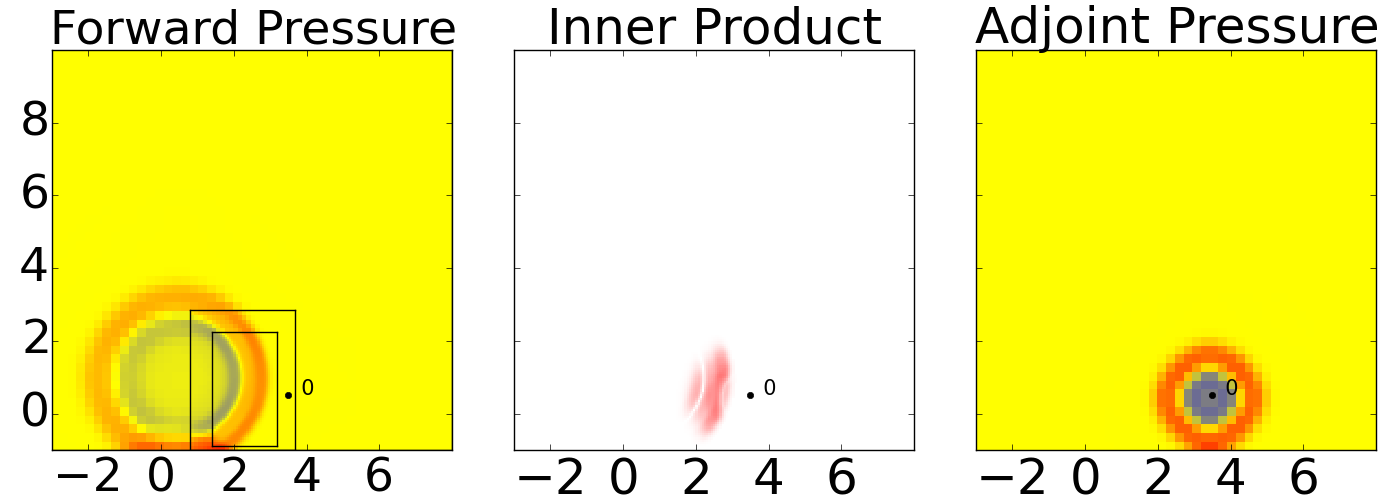}
   \label{fig:2dAcoustics_tp3}
   }
     \subfigure[Results at $t = 1.5$.]{
  \includegraphics[width=0.47\textwidth]{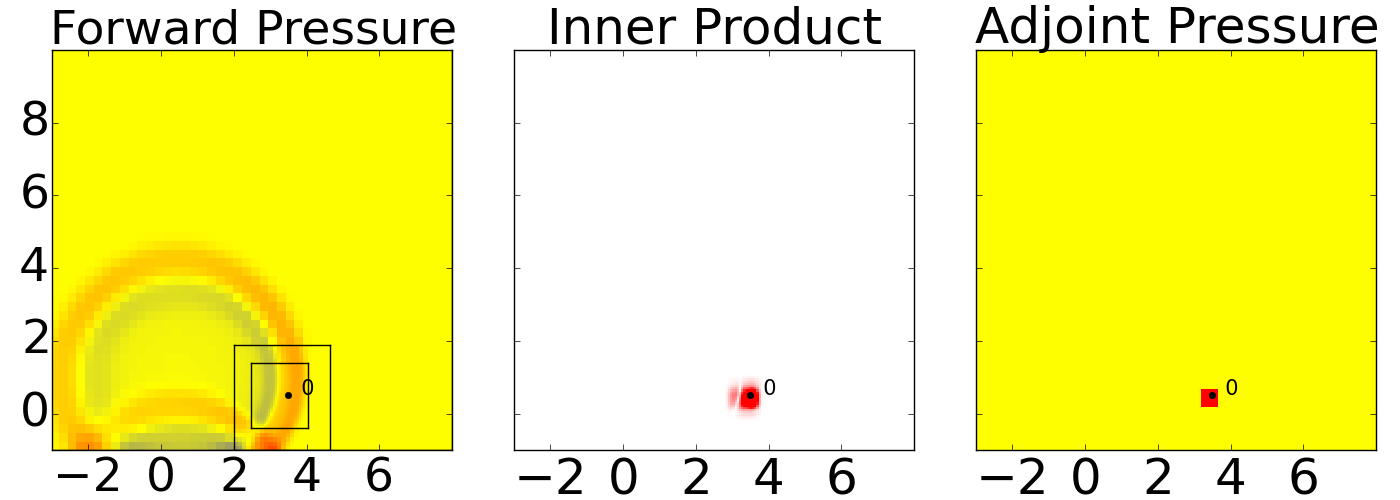}
   \label{fig:2dAcoustics_tp4}
   }
 \caption[Optional caption for list of figures]{%
  Pressure for the two-dimensional acoustics problem with zero boundary 
  conditions when only the pressure
  over the small region shown by the red square in the adjoint pressure 
  plot at the final time t = 1.5 is of interest. The 
  color scale for both pressure figures goes from blue to red, and 
  ranges between $-1.0$ and 1.0 for the forward pressure and between 
  $-0.25$ and 0.25 for the adjoint pressure. The color scale for the inner 
  product figure goes from white to red, and ranges between 0 and 0.25.
  The $y$-axis is the same for all plots in each subfigure.}
   \label{fig:2dAcoustics_tp}
\end{figure}
As the initial data for the adjoint $\hat{q}(x,y,t_f) = \varphi (x,y)$ 
we have a square pulse in pressure, which was described in
Equations \cref{eq:phi_2d} and \cref{eq:delta_2d}.
As time progresses backwards, waves radiates outward and reflect off the
walls.
As the initial data for $q(x,y,t)$ we have 
$q(x,y,0) = \left[ p(x,y,0), 0, 0\right]^T$ 
where $p(x,y,0)$ is given in Equation \cref{eq:2d_acoustics_p}.

\cref{fig:2dAcoustics_tp} shows the Clawpack results for the pressure
for both the forward and adjoint solutions at various times, as well as the 
inner product between the two. Here it is easy to see that the refinement is 
occurring where the inner product is large, and how the interaction between
the forward and adjoint problems at each time step is generating the inner product.

\subsubsection{Wall boundary conditions for a time range of interest}\label{sec:acoustics_reflecting}
Now suppose that $t_f = 6$ seconds, and that we are interested in the accurate 
calculation of the pressure over the time interval $1.0 \leq t \leq 6.0$ 
(so, $t_s = 1$ second). As in the previous example, consider 
wall boundary conditions \cref{eq:acoustics2d_wallsBC}.
Then all but the first term in Equation \cref{eq:acousticseqn_2d} vanish if we 
define the adjoint problem \cref{eq:acoustics2d_adjointwalls}.

\begin{figure}[ht]
 \centering
  \subfigure[$t = 0.0$]{
  \includegraphics[width=0.225\textwidth]{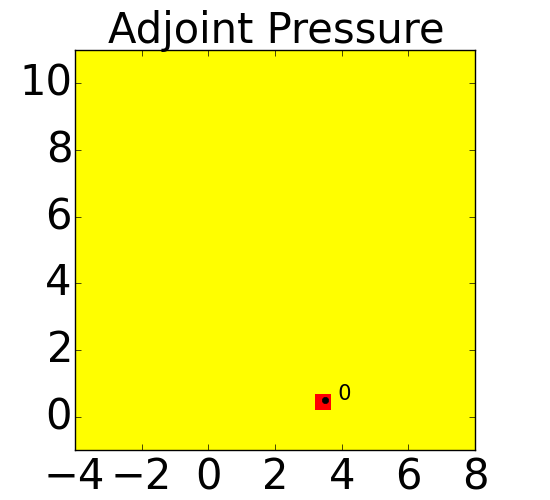}
   \label{fig:2dAcoustics_a1_walls}
   }
   \subfigure[$t = 2.7$]{
  \includegraphics[width=0.225\textwidth]{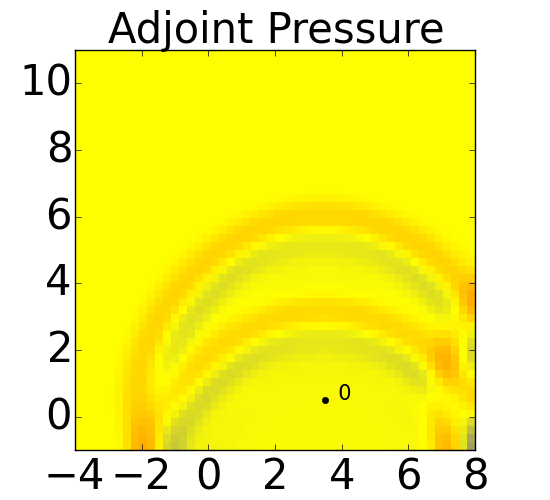}
   \label{fig:2dAcoustics_a2_walls}
   }
  \subfigure[$t = 4.5$]{
  \includegraphics[width=0.225\textwidth]{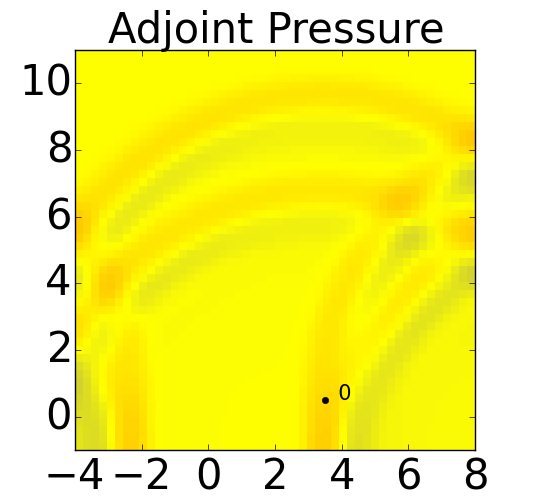}
   \label{fig:2dAcoustics_a3_walls}
   }
  \subfigure[$t = 6.0$]{
  \includegraphics[width=0.225\textwidth]{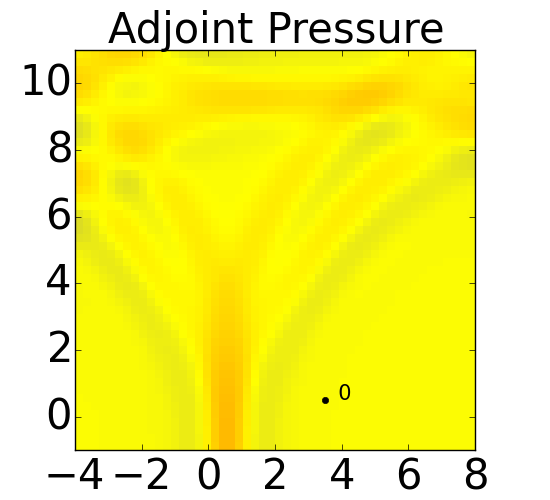}
   \label{fig:2dAcoustics_a4_walls}
   }
 \caption[Optional caption for list of figures]{%
  Pressure for the two-dimensional acoustics adjoint problem with zero 
  boundary conditions. Results are for the shown times given in seconds 
  before the final time, since the ``initial'' conditions are given at the 
  final time.
  The color scale goes from blue to red, and ranges between -0.3 and 0.3.}
   \label{fig:2dAcoustics_a_walls}
\end{figure}

\begin{figure}[ht]
 \centering
  \subfigure[Results at $t = 0.0$ seconds.]{
  \includegraphics[width=0.475\textwidth]{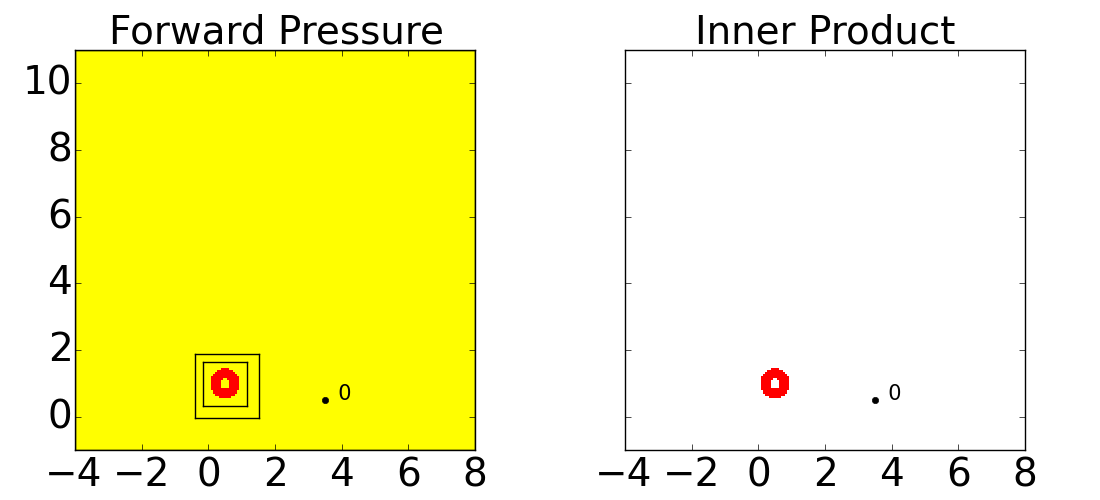}
   \label{fig:2dAcoustics_f1_walls}
   }
   \subfigure[Results at $t = 0.9$ seconds.]{
  \includegraphics[width=0.475\textwidth]{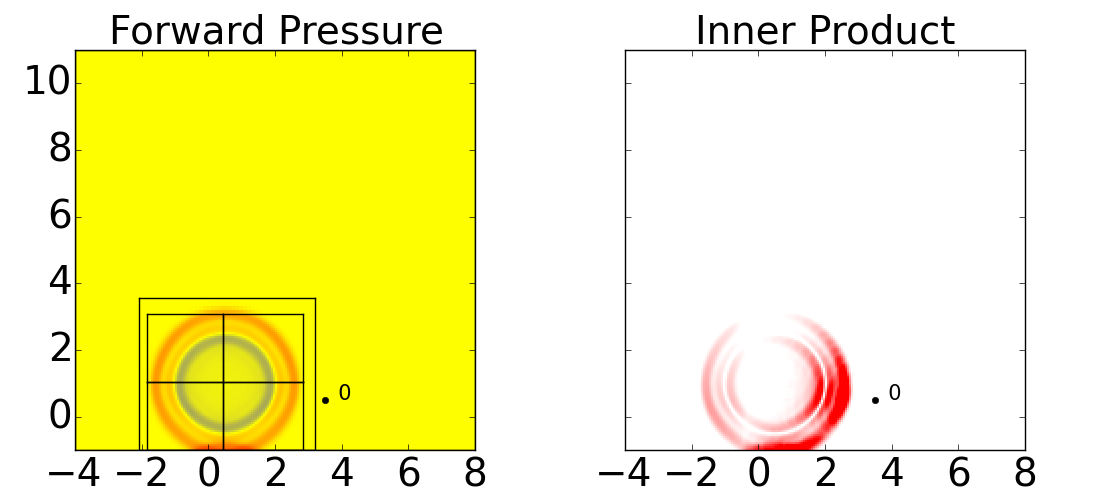}
   \label{fig:2dAcoustics_f2_walls}
   }
  \subfigure[Results at $t = 2.1$ seconds.]{
  \includegraphics[width=0.475\textwidth]{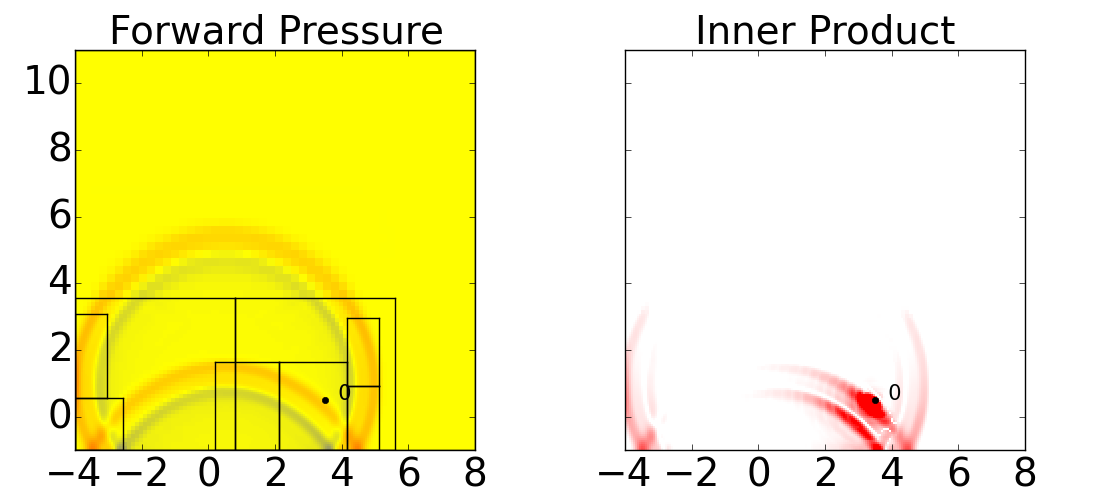}
   \label{fig:2dAcoustics_f3_walls}
   }
  \subfigure[Results at $t = 3.9$ seconds.]{
  \includegraphics[width=0.475\textwidth]{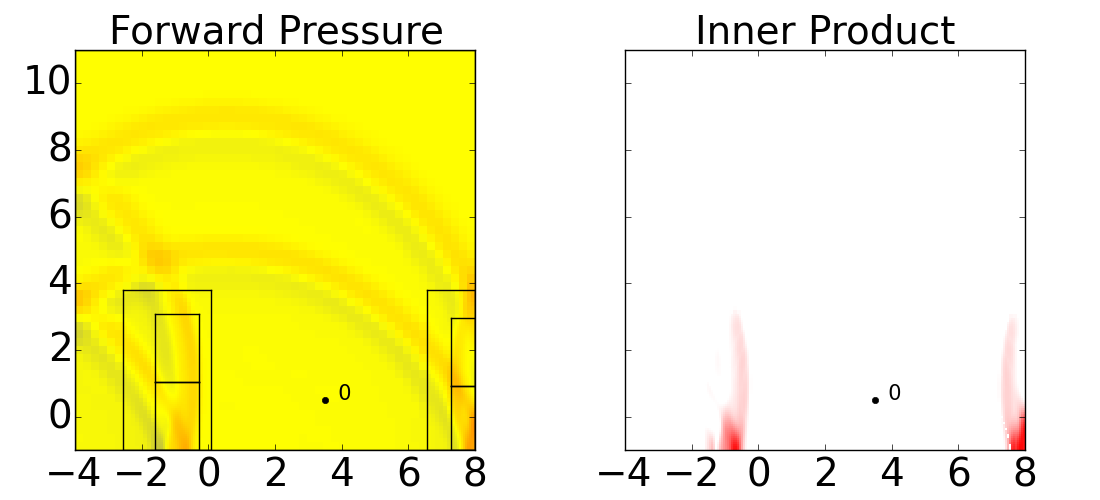}
   \label{fig:2dAcoustics_f4_walls}
   }
     \subfigure[Results at $t = 5.1$ seconds.]{
  \includegraphics[width=0.475\textwidth]{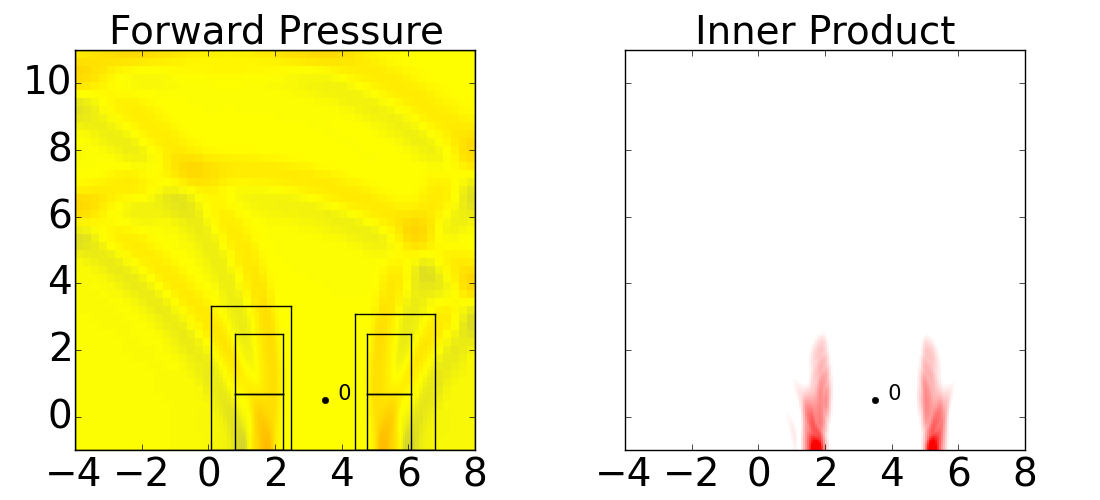}
   \label{fig:2dAcoustics_f5_walls}
   }
     \subfigure[Results at $t = 6.0$ seconds.]{
  \includegraphics[width=0.475\textwidth]{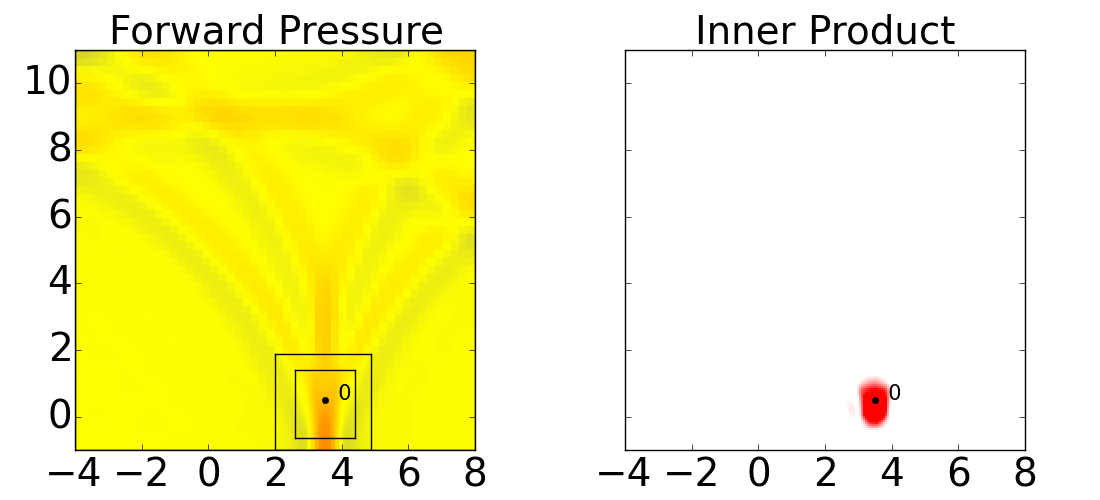}
   \label{fig:2dAcoustics_f6_walls}
   }
 \caption[Optional caption for list of figures]{%
  Computed results for two-dimensional acoustics forward problem 
  with zero boundary conditions. The $y$-axis is the same for both 
  plots in each subfigure. The 
  color scale for the pressure figures goes from blue to red, and 
  ranges between $-1.5$ and 1.5. The color scale for the inner 
  product figures goes from white to red, and ranges between 
  0.01 and 0.12.}
   \label{fig:2dAcoustics_f_walls}
\end{figure}
Using the same initial conditions at in the previous example, 
\cref{fig:2dAcoustics_a_walls} shows the Clawpack 
results for the adjoint pressure, $\hat{p}(x,y,t)$, at
four different times.
Again, as the initial data for $q(x,y,t)$ we have 
$q(x,y,0) = \left[ p(x,y,0), 0, 0\right]^T$ 
where $p(x,y,0)$ is given in Equation \cref{eq:2d_acoustics_p}. To more easily visualize which
areas of the computational domain should be refined at each time step, it is helpful 
to consider the maximum inner product over the appropriate time range,
\begin{align*}
\max\limits_{T \leq \tau \leq t} \hat{q}^T(x,y,\tau ) q(x,y,t),
\end{align*}
with $T = \min (t + t_f - t_s, 0)$. Recall that the areas where this is large are the 
areas where adaptive mesh refinement should take place.  
\cref{fig:2dAcoustics_f_walls} shows the Clawpack results for the pressure, $p(x,y,t)$, 
at six different times, along with this maximum inner product. 

\subsubsection{Mixed boundary conditions for a time range of interest}\label{sec:acoustics_outflow}

Now suppose one of the boundaries has outflow or non-reflecting
boundary conditions, e.g. at $x=a$.  The assumtion in this case is
that this is an artificial computational boundary and ideally we
would be solving a problem with no boundary at this edge.  In other
words we are assuming that if we were to solve the problem on a
larger domain, extended for $x<a$, then waves that pass $x=a$ will
not give rise to reflections that re-enter the computational domain.
If we now consider the adjoint equation on the hypothetical extended
domain, it is clear that waves in the adjoint solution should also
pass through the computational boundary at $x=a$.  Hence non-reflecting
boundary conditions are the correct conditions to impose on the
adjoint solution at any boundary where non-reflecting boundaries
are assumed for the forward problem.

\begin{figure}[ht]
 \centering
  \subfigure[$t = 0.0$]{
  \includegraphics[width=0.225\textwidth]{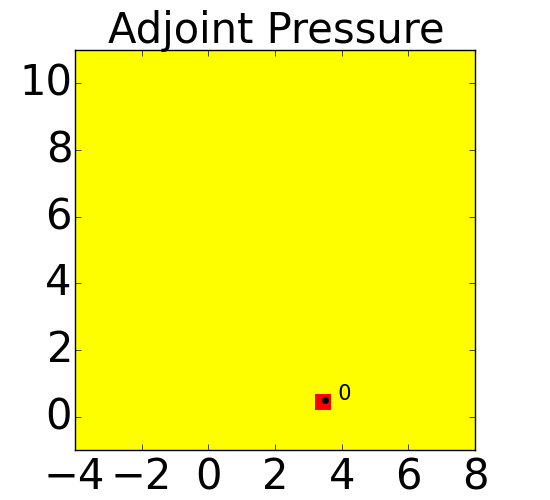}
   \label{fig:2dAcoustics_a1_outflow}
   }
   \subfigure[$t = 2.7$]{
  \includegraphics[width=0.225\textwidth]{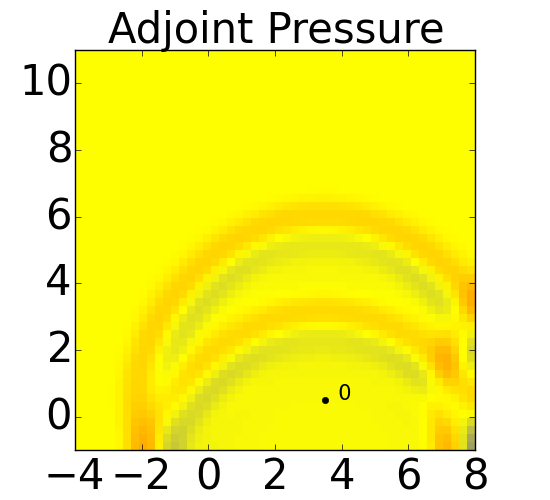}
   \label{fig:2dAcoustics_a2_outflow}
   }
  \subfigure[$t = 4.5$]{
  \includegraphics[width=0.225\textwidth]{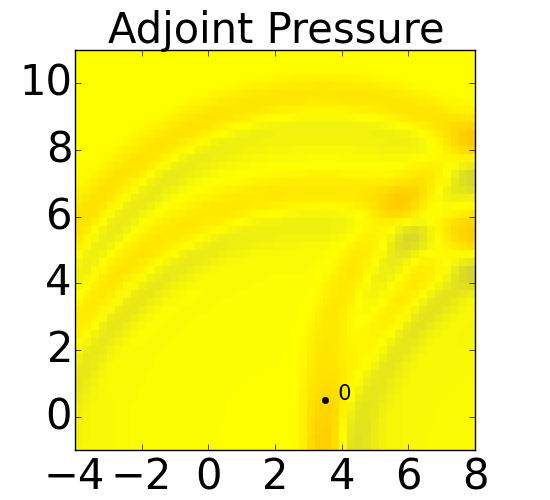}
   \label{fig:2dAcoustics_a3_outflow}
   }
  \subfigure[$t = 6.0$]{
  \includegraphics[width=0.225\textwidth]{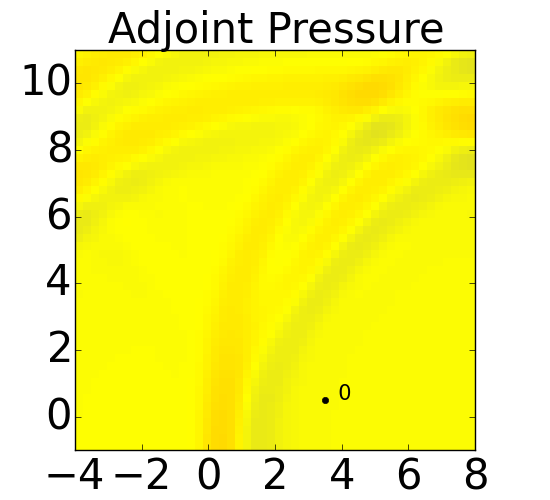}
   \label{fig:2dAcoustics_a4_outflow}
   }
 \caption[Optional caption for list of figures]{%
  Pressure for the two-dimensional acoustics adjoint problem with zero 
  boundary conditions 
  at $x = b$ and $y = \alpha$, and outflow boundary conditions at 
  $x = a$ and $y = \beta$.
  Results are for the shown times given in seconds 
  before the final time, since the ``initial'' conditions are given at the 
  final time.
  The color scale goes from blue to red, and ranges between $-0.3$ and 0.3.}
   \label{fig:2dAcoustics_a_outflow}
\end{figure}

\begin{figure}[ht]
 \centering
  \subfigure[Results at $t = 0.0$ seconds.]{
  \includegraphics[width=0.475\textwidth]{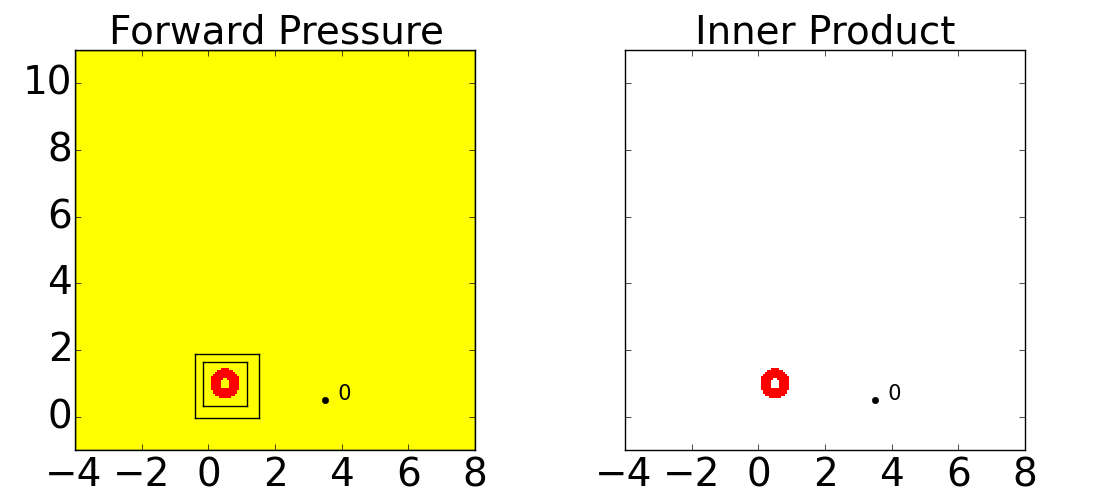}
   \label{fig:2dAcoustics_f1_outflow}
   }
   \subfigure[Results at $t = 0.9$ seconds.]{
  \includegraphics[width=0.475\textwidth]{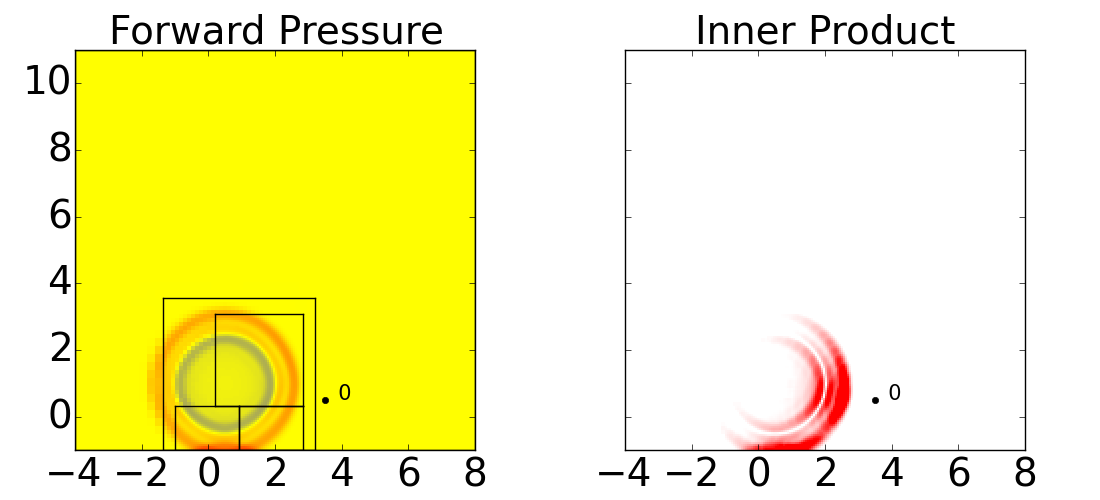}
   \label{fig:2dAcoustics_f2_outflow}
   }
  \subfigure[Results at $t = 2.1$ seconds.]{
  \includegraphics[width=0.475\textwidth]{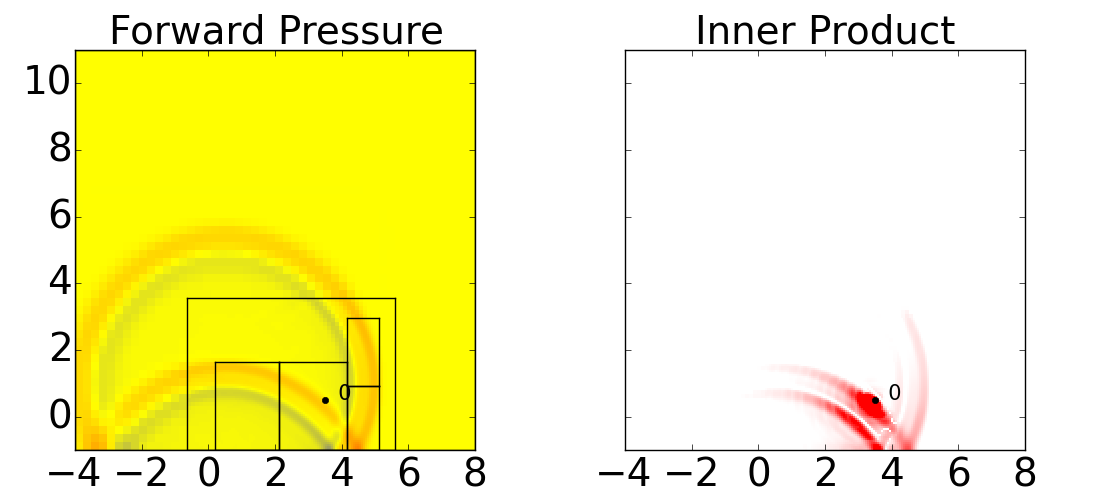}
   \label{fig:2dAcoustics_f3_outflow}
   }
  \subfigure[Results at $t = 3.9$ seconds.]{
  \includegraphics[width=0.475\textwidth]{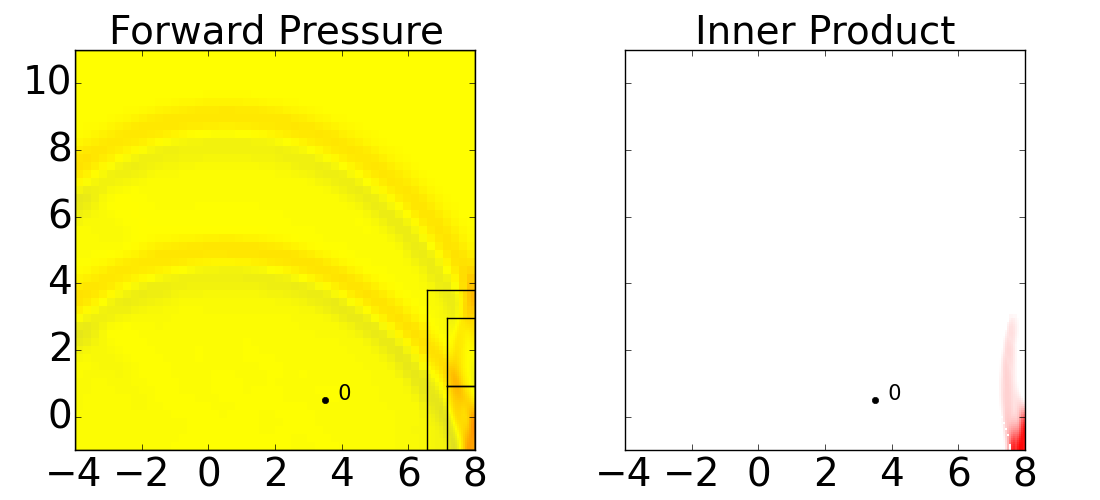}
   \label{fig:2dAcoustics_f4_outflow}
   }
     \subfigure[Results at $t = 5.1$ seconds.]{
  \includegraphics[width=0.475\textwidth]{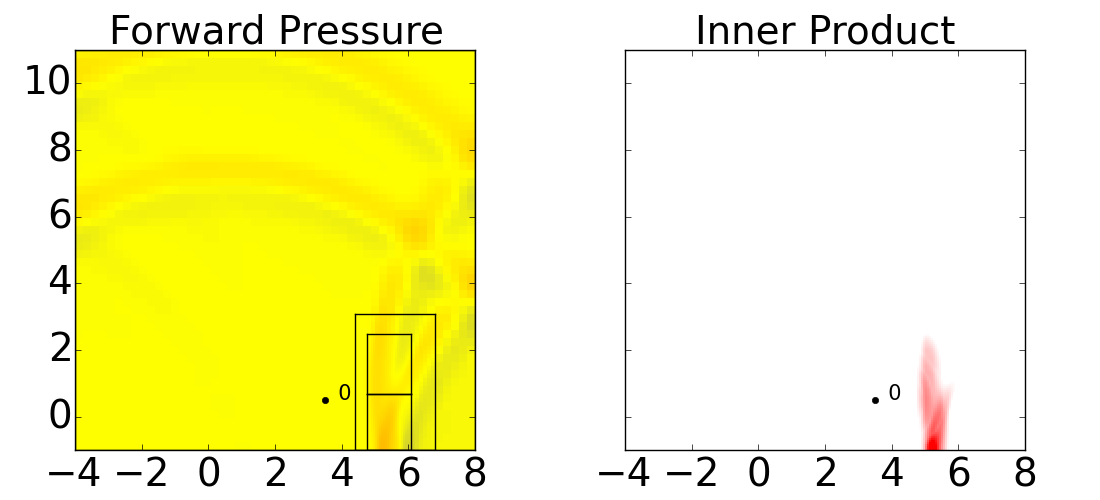}
   \label{fig:2dAcoustics_f5_outflow}
   }
     \subfigure[Results at $t = 6.0$ seconds.]{
  \includegraphics[width=0.475\textwidth]{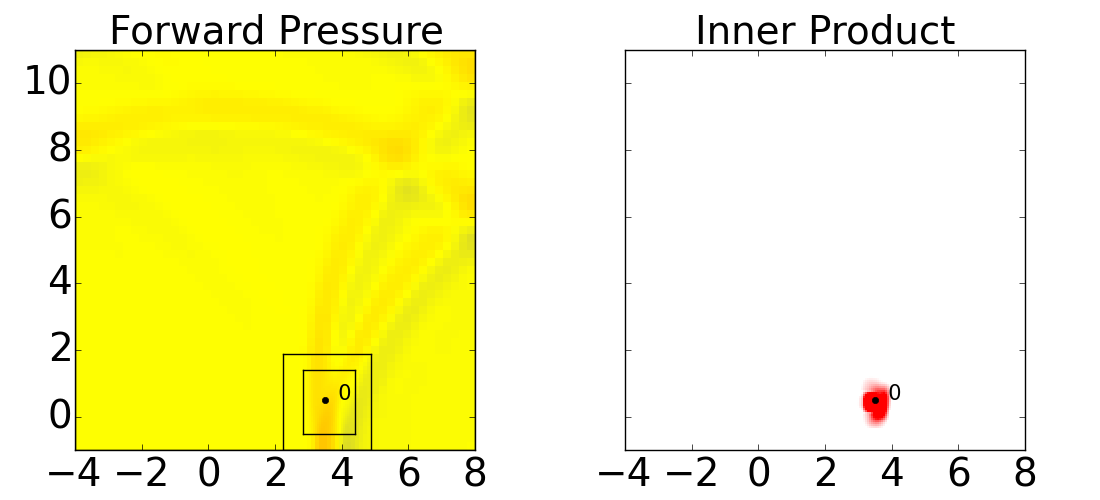}
   \label{fig:2dAcoustics_f6_outflow}
   }
 \caption[Optional caption for list of figures]{%
  Computed results for two-dimensional acoustics forward problem with zero 
  boundary conditions 
    at $x = b$ and $y = \alpha$, and outflow boundary conditions at $x = a$ 
    and $y = \beta$. The $y$-axis is the same for both plots in each subfigure.
    The color scale for the pressure figures goes from blue to red, and 
  ranges between $-1.5$ and 1.5. The color scale for the inner 
  product figures goes from white to red, and ranges between 
  0.01 and 0.12. Compare this to \cref{fig:2dAcoustics_f_walls}, where we 
  were refining around two waves converging on the gauge.}
   \label{fig:2dAcoustics_f_outflow}
\end{figure}

Using the same initial conditions and values for $t_s$ and $t_f$
as in the previous example, we now
replace the boundary conditions at $x=a$ and $y=\beta$ by outflow boundary
conditions to illustrate  this.
\Cref{fig:2dAcoustics_a_outflow} shows the Clawpack 
results for the pressure, $\hat{p}(x,y,t)$, at
four different times.
\cref{fig:2dAcoustics_f_outflow} shows the Clawpack results for the pressure, $p(x,y,t)$, 
at six different times, along with the maximum inner product in the appropriate time range.

\subsubsection{Computational Performance}\label{sec:acou2dperf}
All of the above examples were run on a laptop, and the same examples
utilizing the standard AMR flagging techniques were also run for timing 
comparisons.

For the acoustics examples, the standard refinement technique in AMRClaw 
computes the spatial difference in each direction and for each component of $q$ and 
flags any points where this is greater than some tolerance.  
We will refer to this flagging technique as pressure-flagging. In contrast, for the
adjoint refinement technique we are computing the inner product between $\hat{q}$
and $q$, where the times at which each of these is evaluated are problem specific, 
and flagging any point where this inner product is greater than some tolerance. 
We will refer to this flagging technique as adjoint-flagging. 

\begin{table}[htbp] \label{tb:timing_acoustics}
\caption{Timing comparison for the examples in Section 
\ref{sec:acoustics2d} given in seconds. }
\begin{center}\footnotesize
\renewcommand{\arraystretch}{1.3}
\begin{tabular}{ l  c  | c | c |}
\cline{3-4}
& & \multicolumn{2}{|c|}{\bf Adjoint-Flagging}  \\ \hline
\multicolumn{1}{ | l | }{ \bf Example} & \bf Pressure-Flagging&\bf Forward & \bf Adjoint\\ \hline
\multicolumn{1}{ | l | }{\Cref{sec:acoustics_timepoint}} &3.137 & 1.325 & 0.465 \\
\multicolumn{1}{ | l | }{\Cref{sec:acoustics_reflecting}} & 15.664 & 6.866 & 0.683\\
\multicolumn{1}{ | l | }{\Cref{sec:acoustics_outflow}} & 13.504 & 5.492 & 0.666 \\
\hline
\end{tabular}
\end{center}
\end{table}

For these examples, 
the tolerance for pressure-flagging is set to 0.1 and the tolerance 
for the adjoint-flagging is 
set to 0.02. The timing results for these examples is shown in \cref{tb:timing_acoustics}. 
In the table, the timings in the ``Pressure-Flagging'' column refer to the time required when 
utilizing the pressure-flagging techniques available in AMRClaw, the timings in the 
``Forward'' column refer to the time required to solve the problem when utilizing the 
adjoint-flagging method, and the timings in the ``Adjoint'' column refer to the time required 
to solve the adjoint equations. 
Note that the adjoint-flagging method does require solving two different problems, but 
the total computational time required is still significantly less than the time 
required for the original refinement method.

Another consideration when comparing the adjoint-flagging method with the
pressure-flagging method already in place in AMRClaw is the accuracy of the
results. To test this, a gauge was placed in each example and the output
at that gauge compared across the different methods. 

\begin{figure}[ht]
 \centering
   \subfigure[From example in  
     \cref{sec:acoustics_timepoint}]{
  \includegraphics[width=0.31\textwidth]{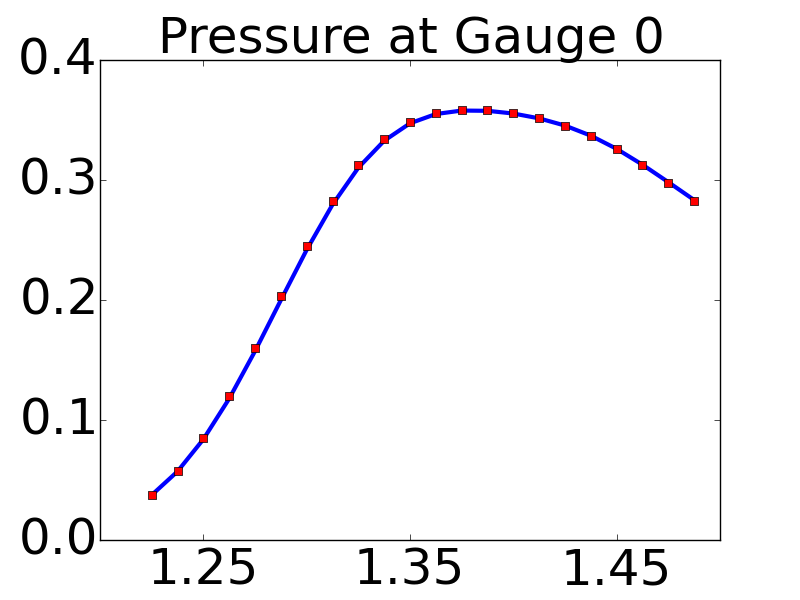}
   \label{fig:gauge_tp}
   }
  \subfigure[From example in 
    \cref{sec:acoustics_reflecting}]{
  \includegraphics[width=0.31\textwidth]{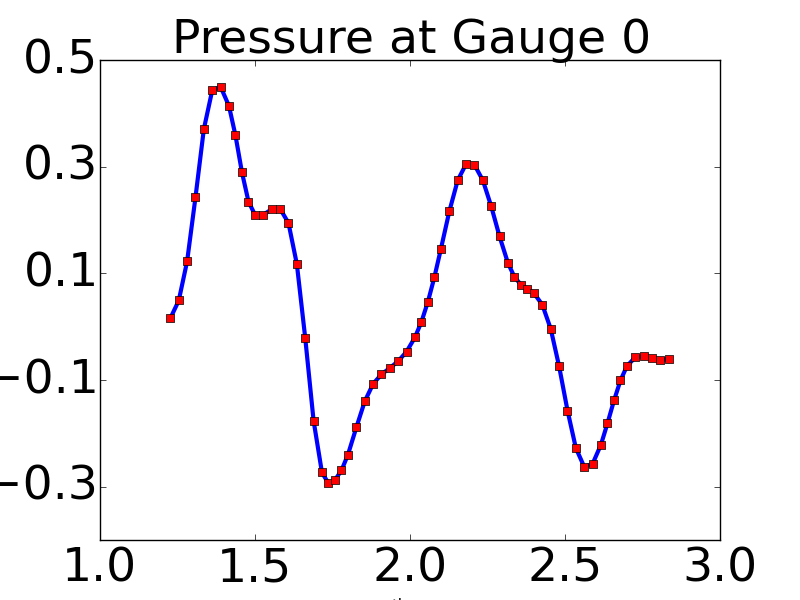}
   \label{fig:gauge_w}
   }
   \subfigure[From example in 
    \cref{sec:acoustics_outflow}]{
  \includegraphics[width=0.31\textwidth]{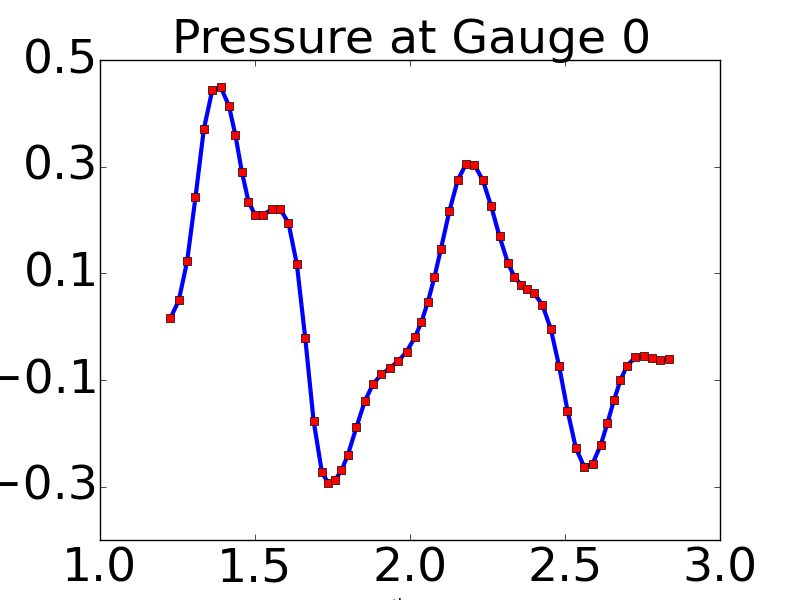}
   \label{fig:gauge_o}
   }
 \caption[Optional caption for list of figures]{%
  Overlaid gauge results for acoustics problems using the
   original and adjoint refinement methods. Results for the original
   pressure-flagging refinement 
   method are shown in red squares, and the results for the adjoint refinement
   method are shown in the blue lines.
   Along the $x$-axis, the time since the beginning of the simulation is shown in seconds. }
   \label{fig:2dAcoustics_gauges}
\end{figure}

In the acoustics examples the gauge was placed at the center of the 
region of interest: $(x,y) = (3.5, 0.5)$, and is shown as a black dot in 
each of the two dimensional acoustics figures previously shown.
 In \cref{fig:2dAcoustics_gauges} the gauge results from the 
adjoint-flagging method are shown as a blue line and then the results from the 
pressure-flagging method were overlaid in red squares. Note that the two 
outputs are indistinguishable for all three examples, indicating that the 
use of the adjoint method did not have a significant detrimental impact on
the accuracy of the calculated results.

\subsection{Tsunami Modeling}\label{sec:tsunami}
Blaise et al. \cite{Blaisea2013} present a method to 
reconstruct the tsunami source using 
a discrete adjoint approach.  In this work we focus on the problem of efficiently
solving the shallow water equations once the initial water displacement is 
known, and discretize the continuous adjoint for the linearized shallow
water equations. The goal is to identify waves propagating
over the ocean that must be accurately tracked because they will reach a
point of interest on the coast, which is otherwise difficult to automate
for reasons mentioned in the introduction.

The algorithms used in GeoClaw for tsunami modeling 
 are described in detail in 
\cite{LeVequeGeorgeBerger2011}, and in general solve the two-dimensional
nonlinear shallow water equations
\begin{align}\label{eq:swe}
h_t + ( hu)_x + ( hv)_y &= 0, \\
( hu )_t + ( hu^2 + \tfrac{1}{2}gh^2)_x + 
( huv)_y &= -ghB_x, \\
( hv)_t + ( huv)_x + ( hv^2 + \tfrac{1}{2}gh^2)_y 
&= -ghB_y,
\end{align}
where $u(x,y,t)$ and $v(x,y,t)$ are the depth-averaged velocities in the two 
horizontal directions, $B(x,y,t)$ is the topography, $g$ is the gravitational 
constant, and $h(x,y,t)$ is the fluid depth. We will use $\eta (x,y,t)$ to denote 
the water surface elevation,
\begin{align*}
\eta (x,y,t) = h(x,y,t) + B(x,y,t),
\end{align*}
$\mu = hu$ to denote the momentum in the $x$ direction, and $\gamma = hv$ 
to denote the momentum in the $y$ direction.

As initial data, we consider a source model for the 1964 Alaska earthquake, 
which generated a tsunami that affected 
the entire Pacific Ocean. The major waves impinging on Crescent City, California from this 
tsunami all occurred within 11 hours after the earthquake, so simulations will be run 
to this time.
 To simulate the effects of this tsunami on Crescent City 
a coarse grid is used over the entire Pacific (1 degree resolution) where the ocean
is at rest. In addition to AMR being used to track propagating waves on finer grids, 
higher levels of refinement are allowed or enforced around Crescent City when the
tsunami arrives. A total of 4 levels of refinement are used, starting with 1-degree 
resolution on the coarsest level, and with refinement ratios of 5, 6, and 6 from one 
level to the next. Only 3 levels were allowed over most of the Pacific, and the remaining
level was used over the region around Crescent City. Level 4, with 20-second resolution,
is still too coarse to provide any real detail on the effect of the tsunami on the harbor. 
It does, however, allow for a comparison of flagging cells for refinement using the adjoint 
method and using the default method implemented in Geoclaw.

\begin{figure}[ht]
 \centering
   \subfigure[Results 2 hours after the earthquake.]{
  \includegraphics[width=0.475\textwidth]{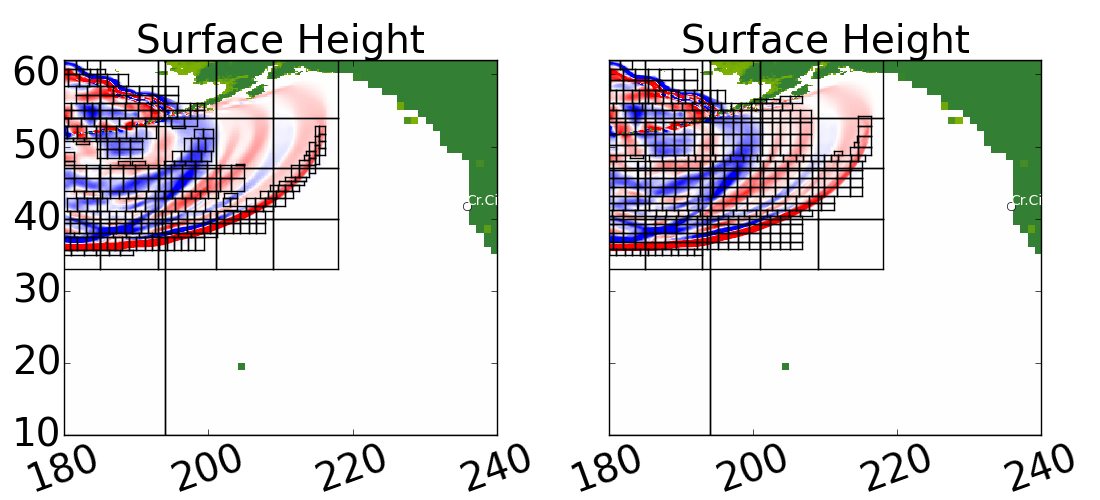}
   \label{fig:Alaska_surfheight_2}
   }
  \subfigure[Results 4 hours after the earthquake.]{
  \includegraphics[width=0.475\textwidth]{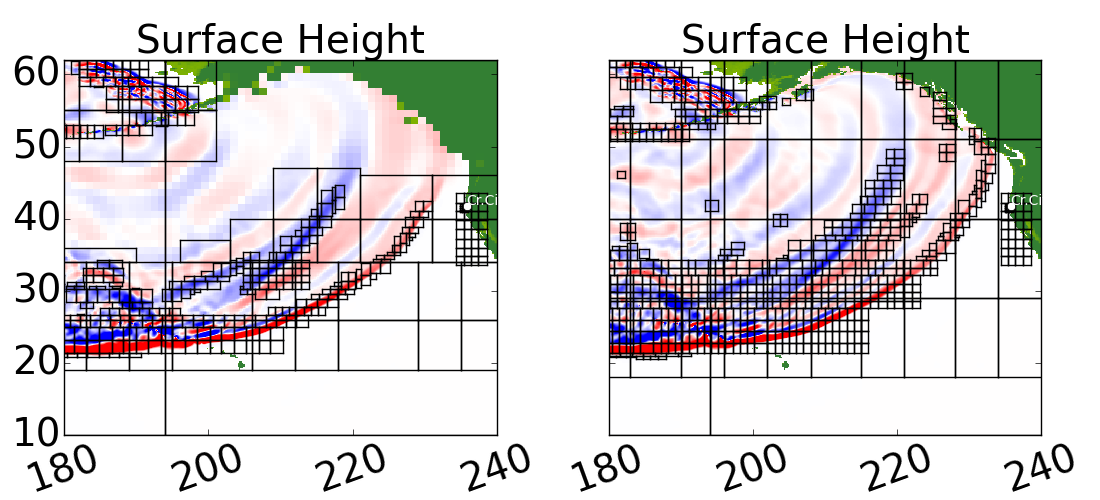}
   \label{fig:Alaska_surfheight_3}
   }
  \subfigure[Results 6 hours after the earthquake.]{
  \includegraphics[width=0.475\textwidth]{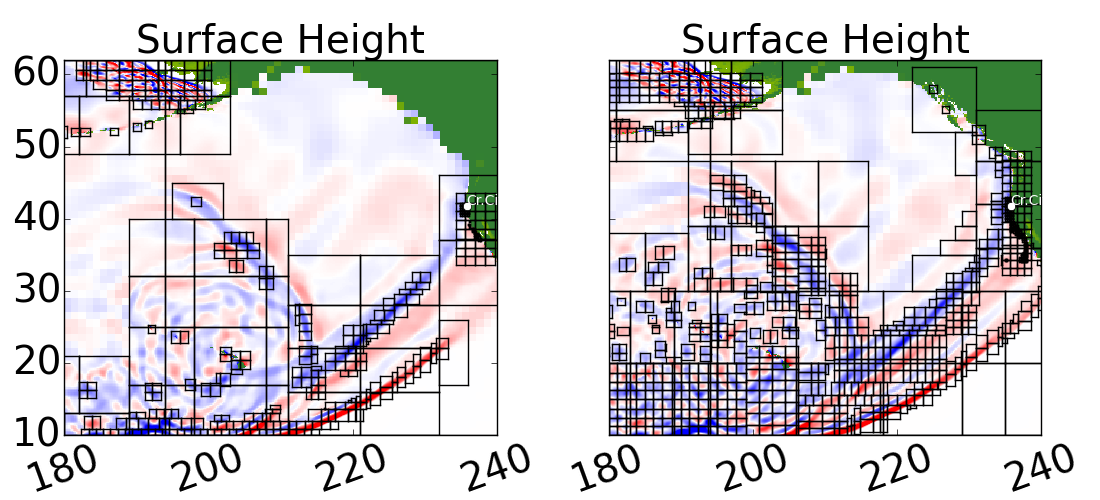}
   \label{fig:Alaska_surfheight_4}
   }
     \subfigure[Results 8 hours after the earthquake.]{
  \includegraphics[width=0.475\textwidth]{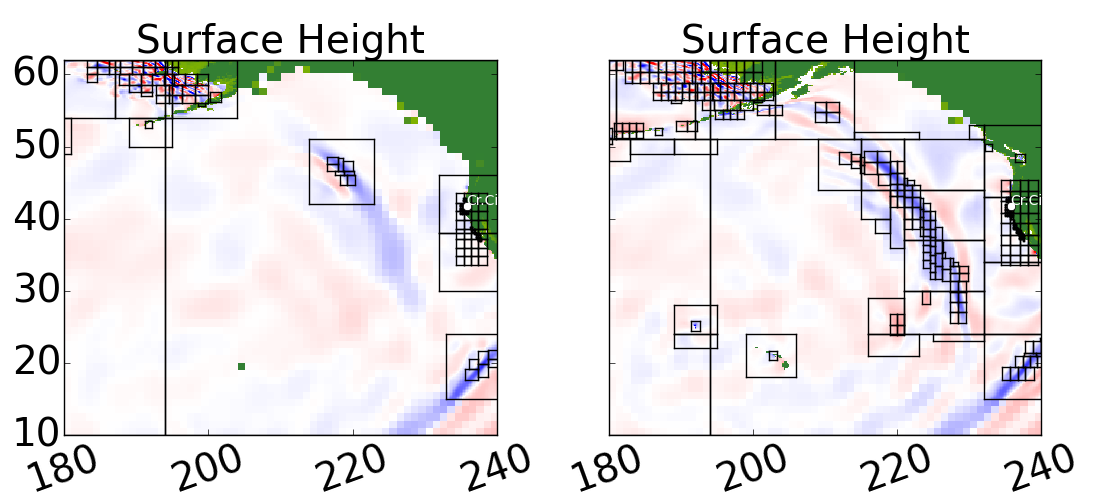}
   \label{fig:Alaska_surfheight_5}
   }
 \caption[Optional caption for list of figures]{%
  Computed results for tsunami propagation problem on two different
  runs utilizing the surface-flagging technique. 
  The $y$-axis is the same for both plots in each subfigure. In each subfigure 
  the left shows results for a tolerance of $0.14$ and the right shows 
  results for a tolerance of $0.09$.
    The color scale goes from blue to red, and 
  ranges between $-0.3$ and 0.3 meters (surface elevation relative to sea
  level).
  Rectangles show patches of up to 60 $\times$ 60 grid cells.}
   \label{fig:Alaska_surfheight}
\end{figure}

The default Geoclaw refinement technique flags cells for refinement when the 
elevation of the sea surface relative to sea level is above some set tolerance 
\cite{LeVequeGeorgeBerger2011}. We will refer to this flagging method as 
surface-flagging. The value selected for this tolerance has a 
significant impact in the results calculated by the simulation, since a smaller tolerance
will result in more cells being flagged for refinement. Consequently, a smaller tolerance
both increases the simulation time required and theoretically increases the accuracy of
the results. 

Two Geoclaw simulations were preformed using surface-flagging, 
one with a tolerance of $0.14$ and another with a tolerance of $0.09$. 
\cref{fig:Alaska_surfheight} shows the results of these two simulations side 
by side for the sake of comparison. Recalling that the black outlines in each figure 
are the edges of grid patches at different
refinement levels, it is easy to see that as expected 
the simulation with the smaller tolerance flagged more cells for refinement. Note that 
the simulation with a surface-flagging tolerance of $0.14$ continues to refine the 
first wave until it arrives
at Crescent City about 5 hours after the earthquake, but after about 6 hours stops 
refining the main secondary wave which reflects off
the Northwestern Hawaiian (Leeward) Island chain before heading 
towards Crescent City. The second simulation, with a surface-flagging tolerance of 
$0.09$ continues to refine this secondary wave until it arrives at Crescent City. 

These two tolerances were selected because they are illustrative of two constraints 
that typically drive a Geoclaw simulation. The larger surface-flagging tolerance, 
of $0.14$, is approximately the largest
tolerance that will refine the initial wave until it reaches Crescent City. Therefore, it 
essentially corresponds to a lower limit on the time required: any Geoclaw simulation with a
larger surface-flagging tolerance would run more quickly but would fail to give accurate results for even
the first wave. Note that for this particular example, even when the simulation will only
give accurate results for the first wave to reach Crescent City, a large area of the wave 
front that is not headed directly towards Crescent City is being refined with the AMR. 
The smaller surface-flagging tolerance, of $0.09$, refines all of the waves of interest that
impinge on Crescent City thereby giving more accurate results at the expense
of longer computational time. 

Now we consider the adjoint approach, which will allow us to refine only those
sections of the wave that will affect our target region.
The linearized shallow water equations govern waves with 
small-amplitude relative to 
the fluid depth. In the case of ocean waves, the nonlinear effects become 
significant when a wave approaches a shore line. Since the water depth in 
the adjoint solution is only used to identify areas of interest, high order 
accuracy is not necessary. Furthermore, computational 
efficiency is the driving force behind our implementation of the adjoint 
method. Therefore, we have solved the simpler linearized equations 
for the adjoint problem. However, the algorithms for the nonlinear 
equations already implemented in GeoClaw are still utilized for the 
forward problem (albeit with a variant of the Roe approximate Riemann 
solver \cite{dgeorge:jcp} that amounts to a local linearization at
each cell edge).

To find the adjoint problem for the linearized equations, we begin by 
letting $(\mu,\gamma) = (hu, hv)$ represent the momentum and
noting that the two momentum equations from \cref{eq:swe} can be rewritten as
\begin{align*}
\mu_t + (hu^2)_x + gh(h+B)_x + (huv)_y &= 0, \\
\gamma_t + (huv)_x + (hv^2)_y + gh(h + B)_y &= 0.
\end{align*}
Linearizing these equations as well as the continuity equation 
about a flat surface $\bar{\eta}$ and zero velocity 
$\bar{u} = 0$, with $\bar{h}(x,y) = \bar{\eta} - B(x,y)$ gives
\begin{align*}
\tilde{\eta}_t + \tilde{\mu}_x + \tilde{\gamma}_y &= 0 \\
\tilde{\mu}_t + g \bar{h}(x,y)\tilde{\eta}_x &= 0 \\
\tilde{\gamma}_t + g\bar{h}(x,y)\tilde{\eta}_y &= 0 
\end{align*}
for the perturbation $(\tilde{\eta}, \tilde{\mu}, \tilde{\gamma})$ about 
$(\bar{\eta}, 0, 0)$. 
Dropping tildes and setting 
\begin{align*}
A(x,y) = \left[ \begin{matrix}
0 & 1 & 0 \\
g \bar{h}(x,y) & 0 & 0 \\
0 & 0 & 0
\end{matrix}\right],\hspace{0.1in}
B(x,y) = \left[ \begin{matrix}
0 & 0 & 1 \\
0 & 0 & 0 \\
g \bar{h}(x,y) & 0 & 0
\end{matrix}\right], \hspace{0.1in}
q(x,y,t) = \left[\begin{matrix}
\eta \\ \mu \\ \gamma
\end{matrix}\right],
\end{align*}
gives us the system $q_t(x,y,t) + A(x,y)q_x(x,y,t) + B(x,y)q_y(x,y,t) = 0$. 

For this example, we are interested in the accurate calculation of the 
water surface height in the area about Crescent City, California. To focus on this
area, we define a circle of radius $1^\circ$
centered about $(x_c,y_c) = (235.80917,41.74111)$ where $x$ and $y$ are being
measured in degrees. Setting 
\begin{align*}
J = 4\int_{x_{min}}^{x_{max}}\int_{y_{min}(x)}^{y_{max}(x)}\eta(x,y,t_f)dy\,dx,
\end{align*}
where the limits of integration define the appropriate circle,
the problem then requires that
\begin{align}\label{eqn:phi_Alaska}
\varphi (x,y) = \left[ \begin{matrix}
I (x,y) \\ 0 \\ 0
\end{matrix}
\right], 
\end{align}
where
\begin{align}\label{eqn:delta_Alaska}
I (x,y) = \left\{
     \begin{array}{ll}
       1 & \hspace{0.3in}\textnormal{if~} \sqrt{\left( x - x_c\right)^2
        + \left( y - y_c\right)^2} \leq 1,\\
       0 & \hspace{0.3in}\textnormal{otherwise.}
     \end{array}
   \right. 
\end{align}
Define
\begin{align*}
\hat{q}(x,y,t_f) = \left[ \begin{matrix}
\hat{\eta}(x,y,t_f) \\ \hat{\mu}(x,y,t_f) \\ \hat{\gamma}(x,y,t_f)
\end{matrix}\right] = \varphi(x,y),
\end{align*}
and note that
\begin{align*}
\int_{t_0}^{t_f}\int_a^b\int_{\alpha}^{\beta}  \hat{q}^T\left( q_t + A(x,y)q_x +
B(x,y)q_y\right) dy\,dx\,dt = 0.
\end{align*}
Again, integrating by parts yields the equation \cref{eq:acousticseqn_2d}
\ignore{
\begin{align}
\int_a^b  \int_{\alpha}^{\beta}  \hat{q}^Tq|^{t_f}_{t_0}&dy\,dx +
\int_{t_0}^{t_f}\int_{\alpha}^{\beta}  \hat{q}^TA(x,y)q|^{b}_{a}dy\,dt 
+ \int_{t_0}^{t_f}\int_{a}^{b}  \hat{q}^TB(x,y)q|^{\beta}_{\alpha}dx\,dt \nonumber \\
&- \int_{t_0}^{t_f}\int_a^b\int_{\alpha}^{\beta}  q^T\left(\hat{q}_{t} +
\left(A^T(x,y)\hat{q}\right)_{x} + \left(B^T(x,y)\hat{q}\right)_{y}\right)
dy\,dx\,dt = 0, \label{eq:shallow}
\end{align}
and if we can define an adjoint problem such that all but the first term 
in this equation vanishes then we are left with
\begin{align*}
\int_a^b \int_{\alpha}^{\beta}  \hat{q}^T(x, y, t_f)q(x, y, t_f) dy\,dx
= \int_a^b \int_{\alpha}^{\beta} \hat{q}^T(x, y, t_0)q(x, y, t_0)dy\,dx,
\end{align*}
which is the expression that allows us to use the inner product of the
adjoint and forward problems at each time step to determine what regions will
influence the point of interest at the final time. 

As we have seen in previous examples, if we define the adjoint equation 
}
and the adjoint equation has the form
\begin{equation}
 \hat{q}_{t} +
\left(A^T(x,y)\hat{q}\right)_{x} + \left(B^T(x,y)\hat{q}\right)_{y} = 0.
\end{equation}
We require this to hold in the same domain as the forward problem, while also
imposing a coastline and boundary conditions similar to those that constrained 
the forward problem.

The topography files used for the adjoint problem are the same as those used
for the forward problem.
However, given that the adjoint problem is being solved on a coarser grid than 
the forward problem, the coastline between the two simulations varies.
Internally, GeoClaw constructs a piecewise-bilinear function 
from the union of any provided topography files. This function 
is then integrated over computational grid cells to obtain a single topography 
value in each grid cell. Note that this means that the topography used in finer
cells for the forward problem has an average value that is equal to the corresponding 
coarse cell value in the adjoint problem. Since the coastline varies between the 
two simulations, when computing the inner product it is possible to find grid cells
that are wet in the forward solution and dry in the adjoint solution. In this case, 
the inner product in those grid cells is set to zero.

The boundary conditions at the coastline also vary slightly between the forward 
and adjoint problems.
The solution of the forward problem involved using a variant of the Roe 
approximate Riemann solver \cite{dgeorge:jcp}, which allows for a changing 
coastline due to inundation. For our simplified solution of the linearized adjoint 
equation, we assume wall boundary conditions (zero normal velocity) at the 
interfaces between any wet cell and dry cell.

\begin{figure}[ht]
 \centering
  \subfigure[1 hour]{
  \includegraphics[width=0.225\textwidth]{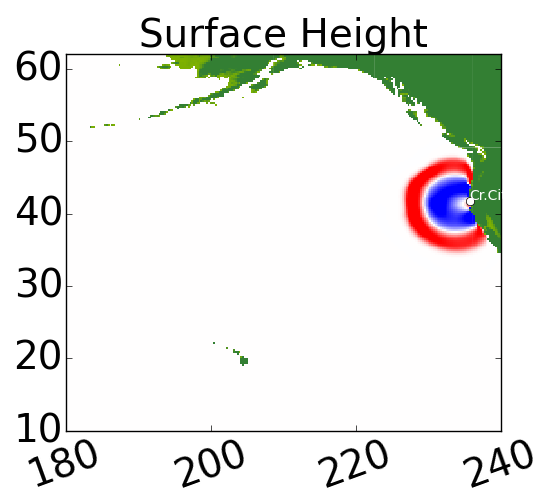}
   \label{fig:Alaska_a1}
   }
   \subfigure[3 hours]{
  \includegraphics[width=0.225\textwidth]{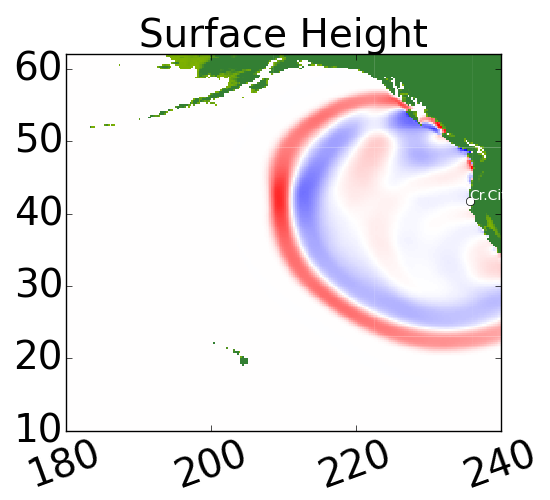}
   \label{fig:Alaska_a2}
   }
  \subfigure[5 hours]{
  \includegraphics[width=0.225\textwidth]{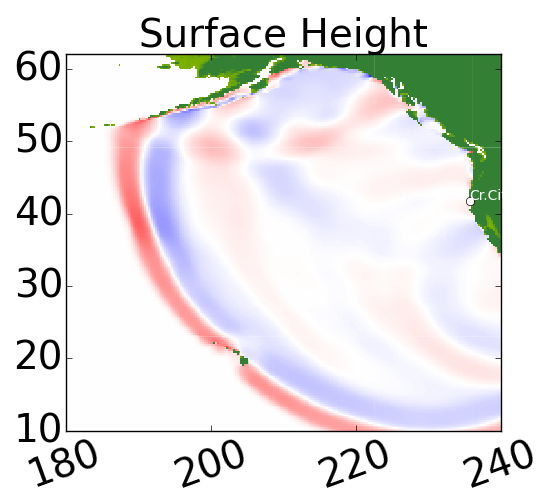}
   \label{fig:Alaska_a3}
   }
  \subfigure[7 hours]{
  \includegraphics[width=0.225\textwidth]{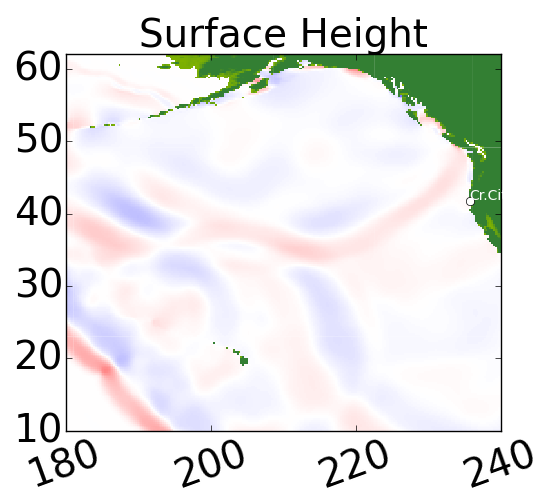}
   \label{fig:Alaska_a4}
   }
 \caption[Optional caption for list of figures]{%
   Computed results for tsunami propagation adjoint problem. Times shown are 
   the number of hours before the final time, since the ``initial'' conditions are 
   given at the final time. 
  The color scale goes from blue to red, and ranges between $-0.05$ and $0.05$.}
   \label{fig:Alaska_adjoint}
\end{figure}

As the initial data for the adjoint $\hat{q}(x,y,t_f) = \varphi  (x,y)$ we have a circular 
area with elevated water height as described in equations \cref{eqn:phi_Alaska} 
and \cref{eqn:delta_Alaska}. 
\cref{fig:Alaska_adjoint} shows the results for the 
simulation of this adjoint problem. For this simulation a grid with 15
arcminute $= 0.25^\circ$
resolution was used over the entire Pacific and no grid refinement was allowed. 
The simulation was run out to 11 hours.
 
The simulation of this tsunami using adjoint-flagging for the AMR
 was run using the same initial grid over the Pacific, the same refinement 
ratios, and the same initial water displacement as our previous 
surface-flagging simulation. The only difference between this simulation
and the previous one is the flagging technique utilized. 
The first waves arrive at Crescent
City around 4 hours after the earthquake, so we set $t_s = 3.5$ hours and $t_f = 11$ 
hours. Recall that the areas where the maximum inner product over the appropriate 
time range,
\begin{align*}
\max\limits_{T \leq \tau \leq t} \hat{q}^T(x,y,\tau ) q(x,y,t)
\end{align*}
with $T = \min (t + t_f - t_s, 0)$, is large are the 
areas where adaptive mesh refinement should take place.

\begin{figure}[ht]
 \centering
  \subfigure[Results 1 hour after the earthquake.]{
  \includegraphics[width=0.475\textwidth]{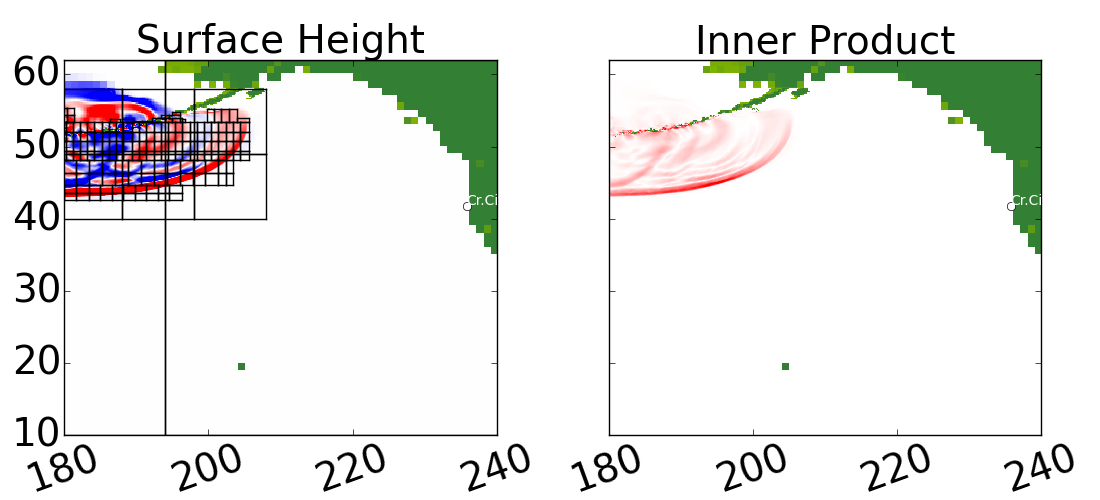}
   \label{fig:Alaska_1}
   }
   \subfigure[Results 2 hours after the earthquake.]{
  \includegraphics[width=0.475\textwidth]{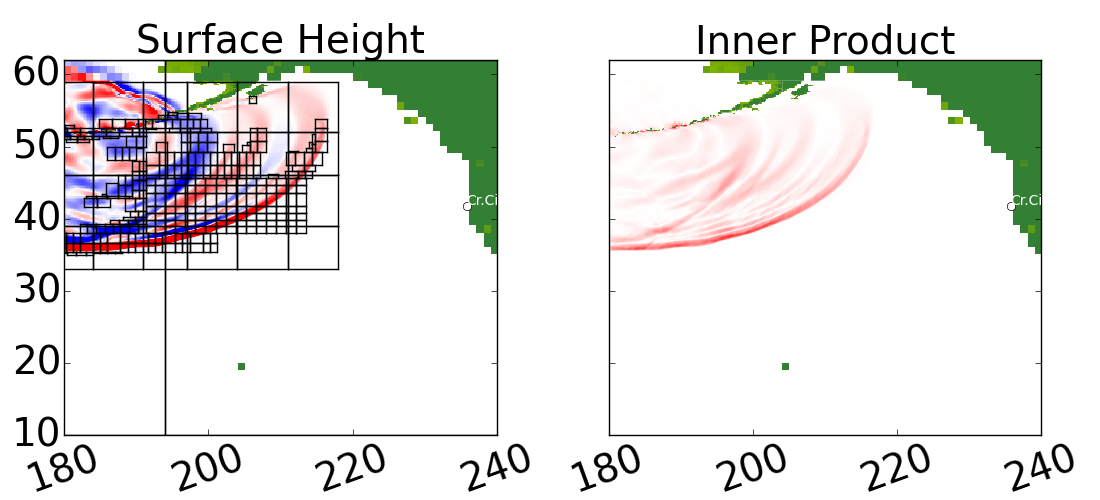}
   \label{fig:Alaska_2}
   }
  \subfigure[Results 4 hours after the earthquake.]{
  \includegraphics[width=0.475\textwidth]{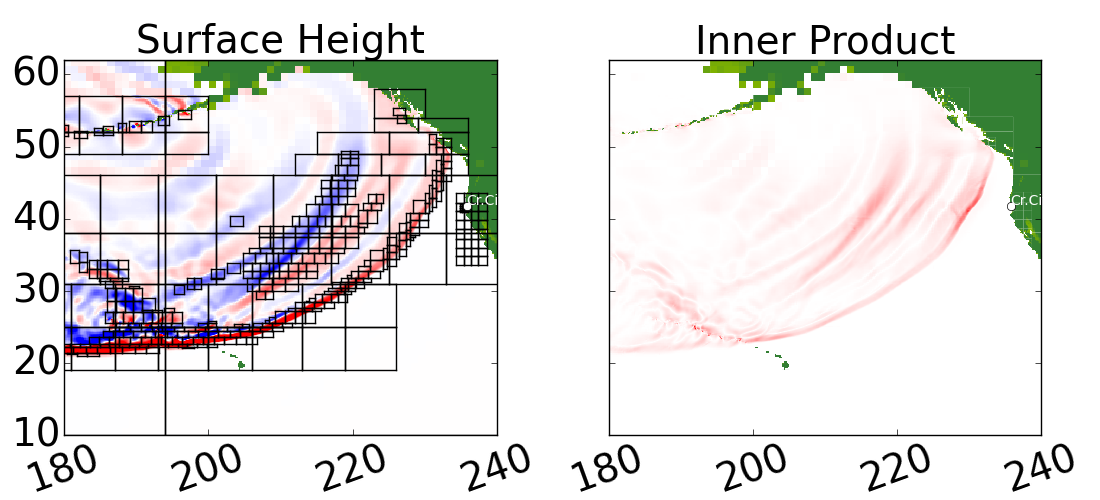}
   \label{fig:Alaska_3}
   }
  \subfigure[Results 6 hours after the earthquake.]{
  \includegraphics[width=0.475\textwidth]{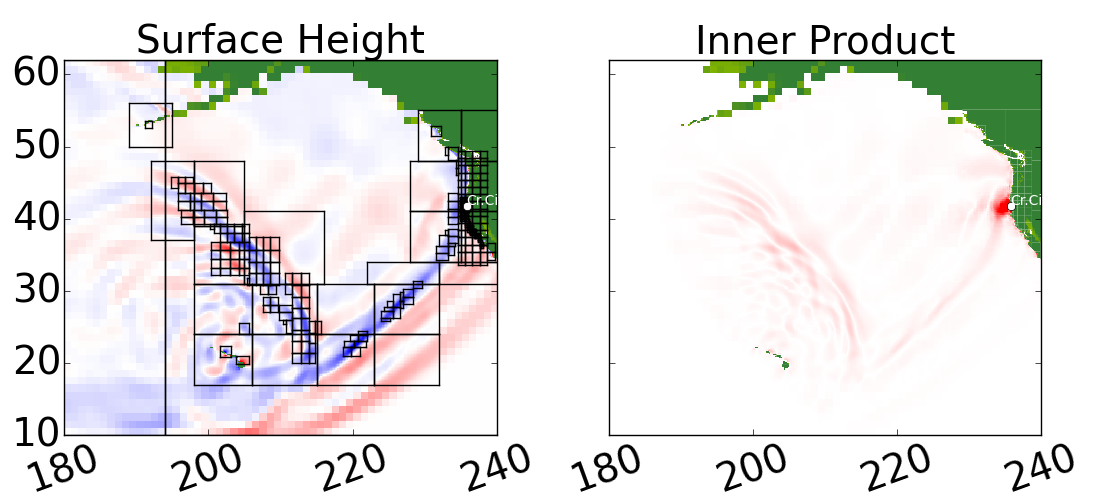}
   \label{fig:Alaska_4}
   }
     \subfigure[Results 8 hours after the earthquake.]{
  \includegraphics[width=0.475\textwidth]{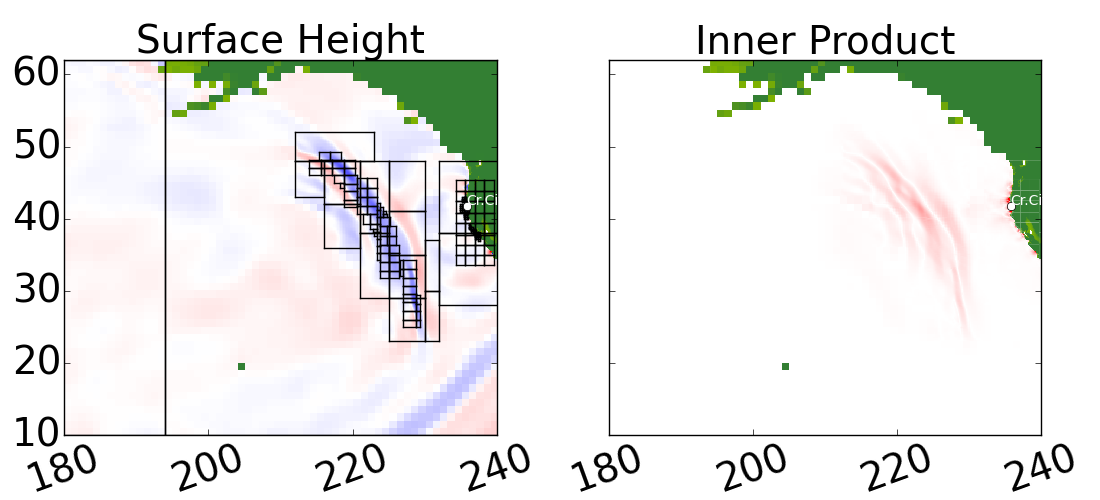}
   \label{fig:Alaska_5}
   }
     \subfigure[Results 10 hours after the earthquake.]{
  \includegraphics[width=0.475\textwidth]{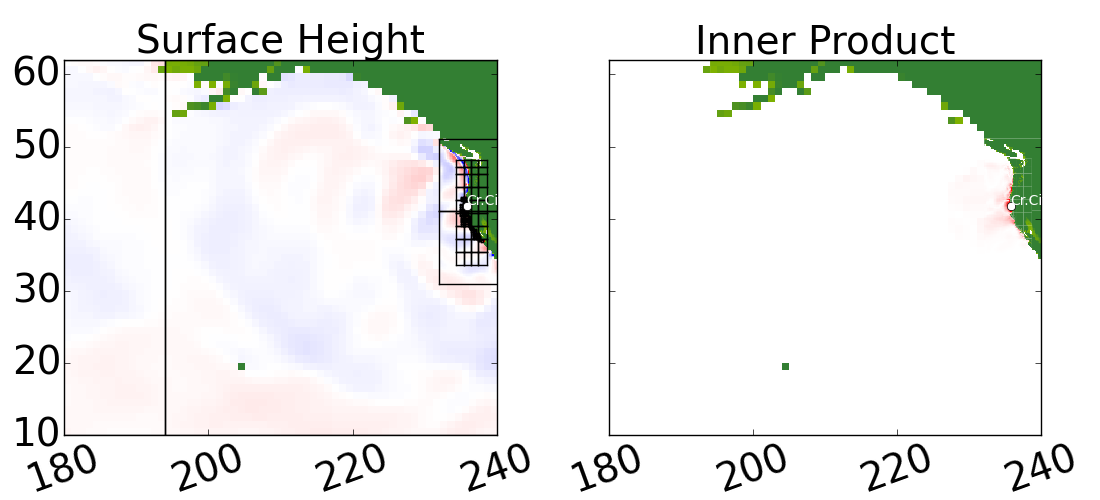}
   \label{fig:Alaska_6}
   }
 \caption[Optional caption for list of figures]{%
  Computed results for tsunami propagation problem. 
  The $y$-axis is the same for both plots in each subfigure.
    The color scale for the surface height figures goes from blue to red, and 
  ranges between $-0.3$ and 0.3. The color scale for the inner 
  product figures goes from white to red, and ranges between 
  0 and 0.05.}
   \label{fig:Alaska}
\end{figure}

\cref{fig:Alaska} shows the Geoclaw results for the surface height at six different times, 
along with the maximum inner product in the appropriate time range. Recalling that the black
outlines in each figure are refined grid patch edges, it is clear that 
this method significantly reduced the extent of the refined regions.

\subsubsection{Computational Performance}\label{sec:tsunamiperf}
The above example was run on a quad-core laptop, for both the surface-flagging 
and adjoint-flagging methods, and the OpenMP option of GeoClaw was enabled 
which allowed all four cores to be utilized.
The timing results for the tsunami simulations are shown in \cref{tb:timing_tsunami}.
 Recall that two simulations were run using 
 surface-flagging, one with a tolerance of $0.14$ (``Large Tolerance'' 
 in the table) and another with a tolerance of $0.09$ (``Small Tolerance'' in the table).
  Finally, a GeoClaw example using adjoint-flagging was run with a 
  tolerance of $0.004$ (``Forward'' in the table), which of course required 
  a simulation of the adjoint problem the timing for which is also shown in the table.
 
 \begin{table}[htbp] \label{tb:timing_tsunami}
\caption{Timing comparison for the example in Section \ref{sec:tsunami} given in seconds. }
\begin{center}\footnotesize
\renewcommand{\arraystretch}{1.3}
\begin{tabular}{| c  | c | c | c|}\hline
\multicolumn{2}{ |c| }{Surface-Flagging} &\multicolumn{2}{ |c| }{Adjoint-Flagging}\\\hline
 \bf Small Tolerance & \bf Large Tolerance &\bf Forward & \bf Adjoint\\ \hline
 8310.285 &  5724.461 & 5984.48 & 26.901\\ \hline
\end{tabular}
\end{center}
\end{table}

  As expected, between the two GeoClaw simulations which utilized  
surface-flagging the one with the larger tolerance took significantly less time. 
Note that although 
  solving the problem using adjoint-flagging did require two different simulations, 
  the adjoint problem and the forward problem, the computational time required falls 
  between the timing required for the two simulations which utilized surface-flagging.

Another consideration when comparing the adjoint-flagging method with the
surface-flagging method already in place in GeoClaw is the accuracy of the
results. To test this, gauges were placed in the example and the output
at the gauges compared across the two different methods. 

\begin{figure}[ht]
 \centering
  \subfigure{
  \includegraphics[width=0.8\textwidth]{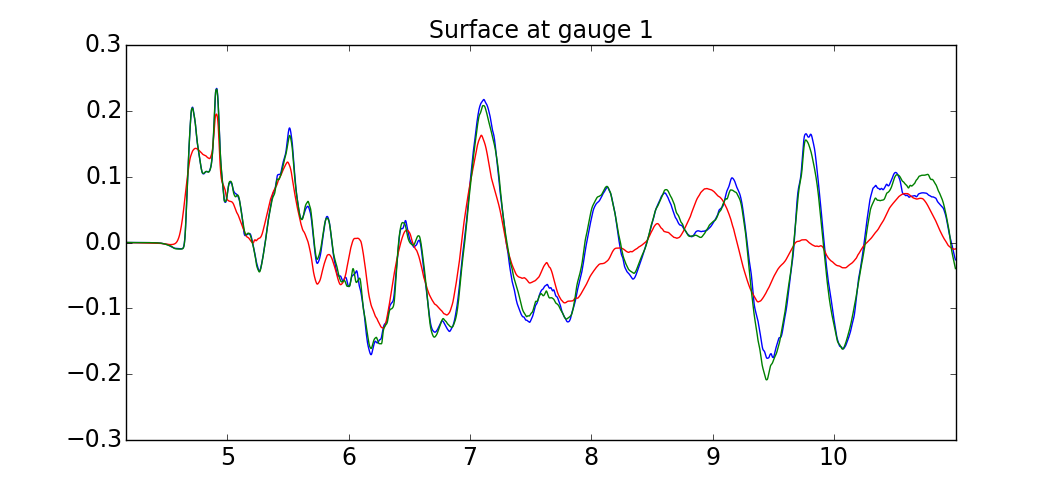}
   \label{fig:Alaska_gauge1}
   }
   \subfigure{
  \includegraphics[width=0.8\textwidth]{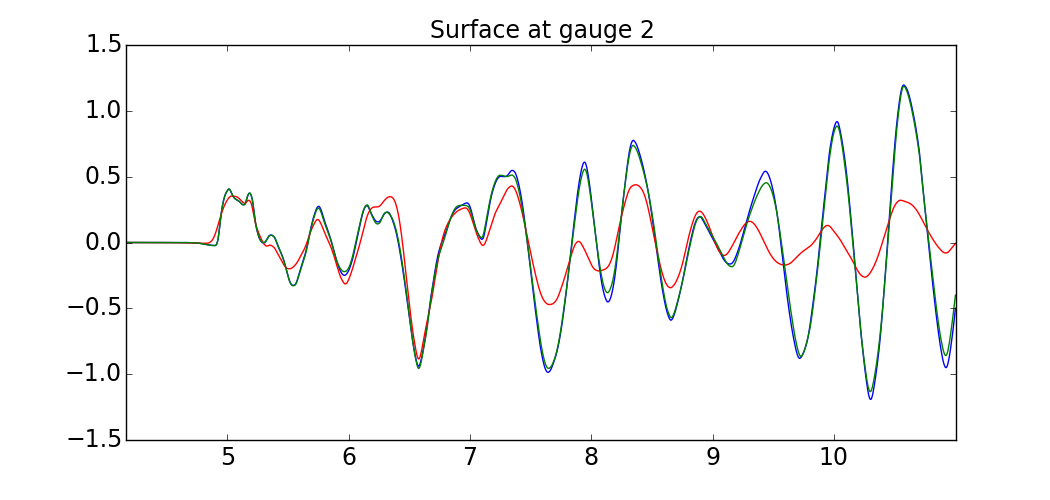}
   \label{fig:Alaska_gauge2}
   }
 \caption[Optional caption for list of figures]{%
  Computed results at gauges for tsunami propagation problem.
  The results from the simulation using the adjoint method are shown in 
  blue, the results from the simulation using the surface-flagging method with 
  a tolerance of $0.14$ are shown in red, and the results from the simulation 
  using the surface-flagging method with a tolerance of $0.09$ are shown in green.
  Along the $x$-axis, the time since the occurrence of the earthquake is shown in hours. 
}
   \label{fig:Alaska_gauges}
\end{figure}

For the tsunami example two gauges are used: gauge 1 is placed at 
$(x,y) = (235.536, 41.67)$ which is on the continental shelf to the west of Crescent 
City, and gauge 2 is placed at $(x,y) = (235.80917,41.74111)$ which is in the 
harbor of Crescent City. In \cref{fig:Alaska_gauges} the gauge results from the adjoint
 method are shown in blue, the results from the surface-flagging technique with a 
tolerance of $0.14$ are shown in red, and the results from the surface-flagging 
technique with a tolerance of $0.09$ are shown in green. Note that the blue and green
lines are in fairly good agreement, indicating that the 
use of the adjoint method did not have a significant detrimental impact on
the accuracy of the calculated results when compared with the smaller tolerance run
using the surface-flagging method. 
However, the larger tolerance run using the surface-flagging method agrees 
fairly well only for the first wave but then rapidly loses accuracy.

\section{Conclusions}\label{sec:conclusions}

Integrating an adjoint method approach to cell flagging into the already existing AMR 
algorithm in Clawpack results in significant time savings for all of the simulations shown 
here. 
For all three two-dimensional acoustics examples examined here, there was no visible loss of accuracy
when using the adjoint method. Furthermore, the pressure-flagging technique required 
 roughly twice as much computational time as the adjoint-flagging technique for each 
 of the examples. It should be noted that these were rather simplified examples chosen 
 for the sake of highlighting different aspects of the adjoint technique. With more complex 
 examples the time savings are likely to be even larger. 
 
For the tsunami example, the surface-flagging technique simulation with a tolerance of
$0.14$ has the advantage of a low computational time but only provides accurate results 
for the first wave to reach Crescent City. The surface-flagging technique simulation
with a tolerance of $0.09$ refines more waves and therefore provides accurate results for 
a longer period of time. However, it has the disadvantage of having a long computational 
time requirement. The use of the adjoint method allows us to retain accurate results 
while also reducing the computational time required.
 
 Therefore, using the adjoint method to guide adaptive mesh refinement can reduce the 
 computational expense of solving a system of equations while retaining the accuracy
 of the results by enabling targeted 
 refinement of the regions of the domain that will influence a specific area of interest.
 
 The code for all the examples presented in this work is available online 
at \cite{adjointCode}, and includes
the code for solving the adjoint Riemann problems. 
This code can be easily modified to solve
different problems. Note that it requires the addition of a new Riemann solver
for the adjoint problem. 
 
 \section{Future Work}\label{sec:futurework}
The method described in this paper flags cells for refinement wherever
the norm of the relevant inner product between the
forward and adjoint solutions is above some tolerance.  Choosing a
sufficiently small tolerance will trigger refinement of
all regions where the forward solution might need to be refined, and in our
current approach these will be refined to the finest level specified in the
computation.  To
optimize efficiency, it would be desirable to have error bounds based on the
adjoint solution that could be used to refine the grid more selectively to
achieve some target error tolerance for the final quantity of interest.
We believe this can be accomplished using  
the Richardson extrapolation error estimator that is built into 
AMRClaw to estimate the point-wise error in the forward solution
and then using the adjoint solution to estimate its effect on the final 
quantity of interest. 
This is currently under investigation and we hope to develop 
a robust strategy that can be applied to a wide variety of problems
for general inclusion into Clawpack.
We will also be
expanding the repository of examples available online.

In this paper the adjoint solution for both the acoustics and tsunami examples was computed on a fixed 
 grid. Allowing AMR to take place when solving the adjoint equation should increase the accuracy 
 of the results, since it would enable a more accurate evaluation of the inner product between the 
 forward and adjoint solutions. 
 In an effort to guide the AMR of the adjoint problem in a similar manner
 to the method used for the forward problem, the two problems would need to be solved somewhat 
 in conjunction and the inner product between the two considered for both the flagging of the 
 cells in the adjoint problem as well as the flagging of cells in the forward problem.
 One approach that is used to tackle this issue is checkpointing schemes, where the forward
 problem is solved and the solution at 
 a small number of time steps is stored for use when solving the 
 adjoint problem \cite{WangMoinIaccarino2009}.
 The automation of this process is another area of future work, 
 and involves developing an evaluation 
 technique for determining the number of checkpoints to save or 
 when to shift from refining the adjoint solution to refining the forward 
 solution and vice versa. 
 
For the tsunami example 
we linearized the shallow water equations about a flat surface, 
which we used to 
solve the adjoint problem. This is sufficient for many important applications, in
particular for tsunami applications where the goal is to track waves in the
ocean where the linearized equtions are very accurate.  If an adjoint
method is desired in the inundation zone, or for other nonlinear hyperbolic
equations,  then the 
adjoint equation is derived by linearizing about
a particular forward solution.  This would 
again require the development of an automated process 
to shift between solving the forward problem, linearizing about that forward problem, and solving 
the corresponding adjoint problem.
Finally, in this work we assumed wall boundary conditions when a wave interacted
with the coastline in the adjoint problem. This assumption, along with the use
 of the linearized shallow water equations, becomes significant when a wave 
 approaches a shore line. Allowing for more accurate iterations between waves 
 and the coastline in the solution of the adjoint problem is another area for 
 future work.

\bibliographystyle{siamplain}
\bibliography{AdjointAMR}

\begin{thebibliography}{10}

\bibitem{AkcelikBirosGhattas2002}
{\sc V.~Akcelik, G.~Biros, and O.~Ghattas}, {\em Parallel multiscale
  {Gauss-Newton-Krylov} methods for inverse wave propagation}, in Proceedings
  of the 2002 ACM/IEEE Conference on Supercomputing, SC '02, 2002, pp.~1--15.

\bibitem{AsnerTavenerKay2012}
{\sc L.~Asner, S.~Tavener, and D.~Kay}, {\em Adjoint-based a posteriori error
  estimation for coupled time-dependent systems}, SIAM J. Sci. Comput., 34
  (2012), p.~A2394–A2419.

\bibitem{BaleLeVequeMitranRossmanith2002}
{\sc D.~S. Bale, R.~J. LeVeque, S.~Mitran, and J.~A. Rossmanith}, {\em A wave
  propagation method for conservation laws and balance laws with spatially
  varying flux functions}, SIAM J. Sci. Comput., 24 (2002), pp.~955--978.

\bibitem{BeckerRannacher2001}
{\sc R.~Becker and R.~Rannacher}, {\em An optimal control approach to a
  posteriori error estimation in finite element methods}, Acta Numerica, 10
  (2001), pp.~1--102.

\bibitem{BergerRigoutsos1991}
{\sc M.~Berger and I.~Rigoutsos}, {\em An algorithm for point clustering and
  grid generation}, IEEE T. Syst. Man Cyb., 21 (1991), pp.~1278--1286.

\bibitem{BergerColella1989}
{\sc M.~J. Berger and P.~Colella}, {\em Local adaptive mesh refinement for
  shock hydrodynamics}, J. Comput. Phys., 82 (1989), pp.~64--84.

\bibitem{BergerGeorgeLeVequeMandli:awr11}
{\sc M.~J. Berger, D.~L. George, R.~J. LeVeque, and K.~T. Mandli}, {\em The
  {GeoClaw} software for depth-averaged flows with adaptive refinement}, Adv.
  Water Res., 34 (2011), pp.~1195--1206,
  \url{www.clawpack.org/links/papers/awr11}.

\bibitem{Berger1998}
{\sc M.~J. Berger and R.~J. LeVeque}, {\em Adaptive mesh refinement using
  wave-propagation algorithms for hyperbolic systems}, SIAM J. Numer. Anal., 35
  (1998), pp.~2298--2316.

\bibitem{BergerOliger1984}
{\sc M.~J. Berger and J.~Oliger}, {\em Adaptive mesh refinement for hyperbolic
  partial differential equations}, J. Comput. Phys., 53 (1984), pp.~484--512.

\bibitem{Blaisea2013}
{\sc S.~Blaise, A.~St-Cyr, D.~Mavriplis, and B.~Lockwood}, {\em {Discontinuous
  Galerkin} unsteady discrete adjoint method for real-time efficient tsunami
  simulations}, J. Comput. Phys., 232 (2013).

\bibitem{BuffoniCupini2001}
{\sc G.~Buffoni and E.~Cupini}, {\em The adjoint advection-diffusion equation
  in stationary and time dependent problems: a reciprocity relation}, Rivista
  di Matematica della Universita di Parma, 4 (2001), pp.~9--19.

\bibitem{BungeHagelbergTravis2003}
{\sc H.-P. Bunge, C.~R. Hagelberg, and B.~J. Travis}, {\em Mantle circulation
  models with variational data assimilation: inferring past mantle flow and
  structure from plate motion histories and seismic tomography}, Geophys. J.
  Int., 152 (2003), pp.~280--301.

\bibitem{Carey2010}
{\sc V.~Carey, D.~Estep, A.~Johansson, M.~Larson, and S.~Tavener}, {\em
  Blockwise adaptivity for time dependent problems based on coarse scale
  adjoint solutions}, SIAM J. Sci. Comput., 32 (2010), pp.~2121--2145.

\bibitem{CLAWPACK}
{\sc {Clawpack Development Team}}, {\em Clawpack software}, 2015,
  \url{http://www.clawpack.org}.
\newblock Version 5.3.

\bibitem{adjointCode}
{\sc B.~N. Davis}, {\em Adjoint code repository},
  \url{https://github.com/BrisaDavis/adjoint}.

\bibitem{dgeorge:jcp}
{\sc D.~L. George}, {\em Augmented {Riemann} solvers for the shallow water
  equations over variable topography with steady states and inundation}, J.
  Comput. Phys., 227 (2008), pp.~3089--3113.

\bibitem{GilesPierce2000}
{\sc M.~B. Giles and N.~A. Pierce}, {\em An introduction to the adjoint
  approach to design}, Flow Turbul. Combust., 65 (2000), pp.~393--415.

\bibitem{Hall1986}
{\sc M.~C.~G. Hall}, {\em Application of adjoint sensitivity theory to an
  atmospheric general circulation model}, J. Atmos. Sci., 43 (1986),
  pp.~2644--2652.

\bibitem{Jameson1988}
{\sc A.~Jameson}, {\em Aerodynamic design via control theory}, J. Sci. Comput.,
  3 (1988), pp.~233--260.

\bibitem{KennedyMartins2013}
{\sc G.~J. Kennedy and J.~R. R.~A. Martins}, {\em An adjoint-based derivative
  evaluation method for time-dependent aeroelastic optimization of flexible
  aircraft}, in Proceedings of the 54th AIAA/ASME/ASCE/AHS/ASC Structures,
  Structural Dynamics, and Materials Conference, Boston, MA, April 2013.
\newblock AIAA-2013-1530.

\bibitem{LeVeque1997}
{\sc R.~J. LeVeque}, {\em Wave propagation algorithms for multidimensional
  hyperbolic systems}, J. Comput. Phys., 131 (1997), pp.~327--353.

\bibitem{Leveque1}
\leavevmode\vrule height 2pt depth -1.6pt width 23pt, {\em Finite Volume
  Methods for Hyperbolic Problems}, Cambridge University Press, 2004.

\bibitem{LeVequeGeorgeBerger2011}
{\sc R.~J. LeVeque, D.~L. George, and M.~J. Berger}, {\em Tsunami modeling with
  adaptively refined finite volume methods}, Acta Numerica,  (2011),
  pp.~211--289.

\bibitem{LiPetzold2004}
{\sc S.~Lia and L.~Petzold}, {\em Adjoint sensitivity analysis for
  time-dependent partial differential equations with adaptive mesh refinement},
  J. Comput. Phys., 198 (2004), pp.~310--325.

\bibitem{MandliDawson2014}
{\sc K.~T. Mandli and C.~N. Dawson}, {\em Adaptive mesh refinement for storm
  surge}, Ocean Modelling, 75 (2014), pp.~36--50,
  \href{http://dx.doi.org/10.1016/j.ocemod.2014.01.002}{doi:10.1016/j.ocemod.2014.01.002},
  \url{http://www.sciencedirect.com/science/article/pii/S1463500314000031}
  (accessed 2015-04-28).

\bibitem{Marburger2012}
{\sc J.~Marburger}, {\em Adjoint-based optimal control of time-dependent free
  boundary problems}, 2012, \url{http://arxiv.org/abs/1212.3789}.

\bibitem{Mishra2013}
{\sc A.~Mishra, K.~Mani, D.~Mavriplis, and J.~Sitaraman}, {\em Time-dependent
  adjoint-based optimization for coupled aeroelastic problems}, in 31st AIAA
  Applied Aerodynamic Conference, San Diego, CA, 2013.
\newblock AIAA 2013-2906.

\bibitem{Nemec}
{\sc M.~Nemec and M.~J. Aftosmis}, {\em Adjoint error estimation and adaptive
  refinement for embedded-boundary {C}artesian meshes}, in 18th {AIAA}
  {C}omputational {F}luid {D}ynamics {C}onference, 2007.

\bibitem{NemecAftosmisWintzer2008}
{\sc M.~Nemec, M.~J. Aftosmis, and M.~Wintzer}, {\em Adjoint-based adaptive
  mesh refinement for complex geometries}, 46th AIAA Aerospace Sciences
  Meeting,  (2008).

\bibitem{Othmer2004}
{\sc C.~Othmer}, {\em Adjoint methods for car aerodynamics}, Journal of
  Mathematics in Industry, 4 (2014), 6.

\bibitem{SandersBradford2002}
{\sc B.~F. Sanders and S.~F. Bradford}, {\em High-resolution, monotone solution
  of the adjoint shallow-water equations}, Int. J. Numer. Meth. Fluids, 32
  (2002), p.~139–161.

\bibitem{SandersKatopodes2000}
{\sc B.~F. Sanders and N.~D. Katopodes}, {\em Adjoint sensitivity analysis for
  shallow-water wave control}, J. Eng. Mech., 126 (2000), p.~909–919.

\bibitem{TrompTapeLie2005}
{\sc J.~Tromp, C.~Tape, and Q.~Liu}, {\em Seismic tomography, adjoint methods,
  time reversal and banana-doughnut kernels}, Geophys. J. Int., 160 (2005),
  pp.~195--216.

\bibitem{VendittiDarmofal2000}
{\sc D.~A. Venditti and D.~L. Darmofal}, {\em Adjoint error estimation and grid
  adaptation for functional outputs: Application to quasi-one-dimensional
  flow}, J. Comput. Phys., 164 (2000), pp.~204--227.

\bibitem{WangMoinIaccarino2009}
{\sc Q.~Wang, P.~Moin, and G.~Iaccarino}, {\em Minimal repetition dynamic
  checkpointing algorithm for unsteady adjoint calculation}, SIAM J. Sci.
  Comput., 31 (2009), pp.~2549--2567.

\end{thebibliography}

\end{document}